\documentstyle[11pt, pb-diagram]{article}
\input amssym.def
\newenvironment{Literatur}[1]{%
\begin{thebibliography}{#1}}{\end{thebibliography}}%

\title{Weak $Spin(9)$-Structures on 16-dimensional Riemannian Manifolds.  \footnote{Supported by the SFB 288 of the DFG.}}
\author{Thomas Friedrich (Berlin)}
\date{\today}

\begin{document}

\newcommand{\D}{\displaystyle}
\newcommand{\upsp}{\phantom{l}}
\newcommand{\downsp}{\phantom{q}}

\maketitle

\mbox{} \hrulefill \mbox{}\\

\newcommand{\vol}{\mbox{vol} \, }
\newcommand{\grad}{\mbox{grad} \, }
\newcommand{\sneun}{Spin(9)}

{\tiny 
\tableofcontents}

\begin{abstract}  The aim of the present paper is the investigation of $\sneun$-structures on 16-dimensional manifolds from the point of view of topology as well as holonomy theory. First we construct several examples. Then we study the necessary topological conditions resulting from the existence of a $\sneun$-reduction of the frame bundle of a 16-dimensional compact manifold (Stiefel-Whitney and Pontrjagin classes). We compute the homotopy groups $\pi_i (X^{84})$ of the space $X^{84}= SO(16) / \sneun$ for $i \le 14$. Next we introduce different geometric types of $\sneun$-structures and derive the corresponding differential equation for the unique self-dual 8-form $\Omega^8$ assigned to any type of $\sneun$-structure. Finally we construct the twistor space of a 16-dimensional manifold with $\sneun$-structure and study the integrability conditions for its universal almost complex structure as well as the structure of the holomorphic normal bundle. 
\end{abstract}

{\small
{\it Subj. Class.:} Differential Geometry.
}\\

\setcounter{section}{0}


\section{Introduction}

The aim of this paper is to present a weak holonomy concept associated to the Lie group $\sneun$. The spin representation of the group $\sneun$ is real and 16-dimensional. According to Berger's holonomy theorem, $\sneun$ can occur as the holonomy group of a 16-dimensional Riemannian manifold. However, D. Alekseevski (see \cite{11x}) and R. Brown/ A. Gray (see \cite{x}) proved that any complete 16-dimensional Riemannian manifold whose holonomy group is contained in $\sneun$ is necessarily flat or isometric to the Cayley plane $F_4 / \sneun$ or its non-compact dual  $F_4^* / \sneun$.\\

In 1971 A. Gray introduced the concept of weak holonomy. He proved that if a manifold has one of the groups
$$ G=SO(n), \quad SU(n), \quad Sp(n) \cdot Sp(1), \quad Sp(n) \cdot SO(2) , \quad Sp(n), \quad Spin(7) $$

as weak holonomy group, then its holonomy is in fact already contained in $G$. Consequently, only 3 groups may admit a weak holonomy concept that is more general than the traditional holonomy approach (see \cite{Gray2}):
$$ G= U(n), \quad G=G_2 \, \, \mbox{in dimension 7} , \quad G= \sneun \, \, \mbox{in dimension 16} . $$

The first two cases yield a rich geometric structure both as weak and as classical holonomy groups and have been studied intensively. Manifolds with weak holonomy group $U(n)$ are called nearly K\"ahler (see \cite{Gray1}). A. Gray has investigated them since 1976 (see \cite{Gray3}) and pointed out that they have special properties in dimension 6. This effect is closely related to the fact that, on a 6-dimensional manifold, the existence of a nearly K\"ahler structure is equivalent to the existence of a real Killing spinor (see \cite{Grune}). In 1981 S. Marchiafava (see \cite{March}) characterized 7-dimensional manifolds with weak holonomy group $G_2$, and M. Fernandez / A. Gray (see \cite{FG1}) studied the different geometric types of $G_2$-structures systematically. In particular, nearly parallel $G_2$-structures correspond  again to real Killing spinors (see \cite{FKMS}). Only the case of weak holonomy $\sneun$ on 16-dimensional manifolds has been neglected until now. \\

We will first  define a (topological) $\sneun$-structure on a 16-dimensional manifold as some 9-dimensional subbundle $V^9$ of the bundle of endomorphisms $\mbox{End} (T(M^{16}))$. Locally, there exist nine endomorphisms $I_{\alpha} \in \Gamma (V^9) \, \, (1 \le \alpha \le 9)$ satisfying the relations
$$ I_{\alpha}^2 = Id \quad , \quad I_{\alpha}^* = I_{\alpha} \quad , \quad I_{\alpha} I_{\beta} = - I_{\beta} I_{\alpha} \quad \mbox{for $ \alpha \not= \beta$} . $$

From this point of view a $\sneun$-structure is a 16-dimensional analogue of a quaternionic structure. It was already noticed in \cite{x} that there exists a $\sneun$-invariant and self-dual 8-form $\Omega^8$ on ${\Bbb R}^{16}$. Although it cannot be used to uniquely characterize the structures we are interested in, it will play an important role. We construct several examples of 16-dimensional manifolds admitting natural topological $\sneun$-structures. Then we derive necessary conditions for the Stiefel-Whitney and the Pontrjagin classes of a compact manifold admitting a $\sneun$-reduction of the frame bundle. For example, the complete intersection of three quadrics in ${\Bbb P}^{11} ({\Bbb C})$ satisfies all these conditions. Up to now it seems impossible to formulate a necessary and sufficient criterion for the existence of a $\sneun$-structure. This is mainly due to the complicated homotopy type of the space $X^{84} = SO(16)/ \sneun$. Using recent results on the homotopy groups $\pi_i (SO(n))$ outside the stable range (see \cite{9x}) we compute $\pi_i (X^{84})$ for $i= 1, \ldots , 14$.\\

In Section 8 we start with the investigation of the geometry of $\sneun$-structures. For this, we assign to any $\sneun$-reduction a 1-form $\Gamma$ with values in the bundle $\Lambda^3 (V^9)$, i.e.,
$$ \Gamma \in \Lambda^1 (M^{16}) \otimes \Lambda^3 (V^9)  . $$

The space $\Lambda^1 ({\Bbb R}^{16}) \otimes \Lambda^3 ({\Bbb R}^9)$ decomposes under the action of $\sneun$ into 4 irreducible summands. Depending on the algebraic type of $\Gamma$ there are 16 different geometric types of $\sneun$-structures. One of the components in the splitting of $\Lambda^1 ({\Bbb R}^{16}) \otimes \Lambda^3 ({\Bbb R}^9)$ is the representation $\Lambda^1 ({\Bbb R}^{16})$ itself. We call the corresponding type of $\sneun$-structure nearly parallel, i.e., a topological $\sneun$-structure is nearly parallel if and only if $\Gamma$ is a vector field. We prove  that $S^1 \times S^{15}$ admits a nearly parallel $\sneun$-structure, thus showing that  such structures do exist. \\

Using the fact that the $\sneun$-representation $\Lambda^7 ({\Bbb R}^{16})$ is multiplicity-free (see \cite{agricola}), we can prove that the above introduced 8-form $\Omega^8$ of a nearly parallel $\sneun$-structure satisfies the equations
$$ \delta \Omega^8 = - 504 (\Gamma \, \mbox{}_{\_\!\!\_\!}{\scriptstyle{\mid}} \, \Omega^8) \quad , \quad d \Omega^8 = - 504 \star (\Gamma \, \mbox{}_{\_\!\!\_\!}{\scriptstyle{\mid}} \, \Omega^8) . $$

The other geometric types of $Spin(9)$-structures yield similar differential equations for $d \Omega^8$.\\

In the final part of this paper we sketch the twistor theory for nearly parallel $\sneun$-structures. For this, we introduce the space ${\cal T}_1$ of all complex structures in $\Lambda^2 ({\Bbb R}^9)$
$$ {\cal T}_1 = \{ {\cal J} = \sum\limits^9_{1 \le \alpha < \beta \le 9} \mbox{x}_{\alpha \beta} I_{\alpha} I_{\beta} : \quad {\cal J}^2 = - Id \} , $$

which is isomorphic to a complex quadric $Q$ in ${\Bbb P}^8 ({\Bbb C})$ and on which $\sneun$ acts transitively. Using ${\cal T}_1$ as a typical fibre we define a twistor space ${\cal T}_1 (M^{16})$ for any 16-dimensional manifold with a fixed $\sneun$-structure. It has a canonical almost complex structure as well as an anti-holomorphic involution without fixed points. From the general theory of twistor spaces (see {\cite{BBO}, \cite{OBR}) we know that the almost complex structure on ${\cal T}_1 (M^{16})$ has to satisfy two integrability conditions. The first one concerns the torsion tensor and turns out to be automatically satisfied in case of a nearly parallel $\sneun$-structure. Therefore, only the integrability condition involving the Riemannian curvature $\Omega^Z$ and the derivative $D^Z (\Gamma)$ has to be fulfilled. As an example, we show that the twistor space of the Cayley plane $F_4 / \sneun$ is isomorphic to
$$ {\cal T}_1 (F_4 / \sneun)=F_4 / (Spin(2) \times_{{\Bbb Z}_2} Spin(7))  . $$

Since $Spin(2) \times_{{\Bbb Z}_2} Spin(7)$ is the centralizer of the subgroup $Spin(2)$ in $F_4$, ${\cal T}_1 (F_4 / \sneun)$ is a generalized flag manifold and therefore a complex projective variety (see \cite{wallach}). The twistor space ${\cal T}_1 (S^1 \times S^{15})$ of $S^1 \times S^{15}$ with its invariant nearly parallel $Spin(9)$-structure is a complex subvariety of the twistor space of $S^1 \times S^{15}$ considered as a conformally flat 16-dimensional Riemannian manifold. Then we  describe the twistor space of the flat manifold ${\Bbb R}^{16}$ as an 8-dimensional holomorphic vector bundle $N$ over the quadric $Q$, and compute its Chern classes as well as the space of all holomorphic sections ${\cal H}^0 (N)$ of this bundle. This result allows not only the description of ${\cal T}_1 ({\Bbb R}^{16})$, but also of the normal bundle $N$ to any fibre inside an arbitrary twistor space ${\cal T}_1 (M^{16})$. It turns out that $N$ admits a 16-dimensional family of holomorphic sections, i.e., it is possible to reconstruct the given manifold $M^{16}$ with a nearly parallel $\sneun$-structure from its twistor space ${\cal T}_1 (M^{16})$.\\

The author thanks Ilka Agricola for her helpful comments and Heike Pahlisch for her competent and efficient {\LaTeX} work.\\

\section{$Spin(9)$-structures on 16-dimensional manifolds}

\newcommand{\rneun}{{\Bbb R}^9}
\newcommand{\cneun}{{\cal C}_9}
\newcommand{\dneun}{\Delta_9}
\newcommand{\kneun}{\kappa_9}

\newfont{\graf}{eufm10}
\newcommand{\g}{\mbox{\graf g}}
\newcommand{\h}{\mbox{\graf h}}
\newcommand{\m}{\mbox{\graf m}}
\newcommand{\n}{\mbox{\graf n}}
\newcommand{\so}{\mbox{\graf so}(16)}

Let $\rneun$ be the 9-dimensional Euclidean space and denote by ${\cal C}_9$  the real Clifford algebra of the negative definite quadratic form. $\cneun$ is generated by the vectors of ${\Bbb R}^9$, and the relation
$$ v \cdot w + w \cdot v = - 2 \langle v,w \rangle \quad , \quad v,w \in \rneun  , $$

holds. The spin representation $\kappa_9$ of the group $\sneun$ is a faithful real representation in the 16-dimensional space $\dneun$ of real spinors and is the unique irreducible representation of the group $\sneun$ in dimension 16. Moreover, $\sneun$ acts transitively on the 15-dimensional sphere $S(\dneun)$ of all spinors of length one. The representation $\kappa_9$ is the isotropy representation of the Cayley plane $F_4 / \sneun$, the unique exceptional symmetric space of rank one. Denote by $\lambda_{16} : {Spin}(16) \to SO (16) =SO(\dneun)$ the universal covering of the orthogonal group $SO (16)$.  $\sneun$ is a simply connected group and there exists a lift $\widetilde{\sneun} \subset Spin(16)$ of the group $\kappa_9 (\sneun)$. Since the subgroup $\sneun \subset SO(16)$ contains the element $(- \mbox{Id})$, the subgroup $\widetilde{\sneun} \subset Spin (16) \subset {\cal C}_{16}$ has to contain one and only one of the two elements $\pm e_1 \cdot \ldots e_{16}$ of the Clifford algebra ${\cal C}_{16}$ of ${\Bbb R}^{16}$. \\

{\bf Proposition 1:} {\it The subgroup $\widetilde{\sneun} \subset Spin(16)$ contains the element $e_1 \cdot \ldots \cdot e_{16}$, but it does not contain the element $- e_1 \cdot \ldots \cdot e_{16}$.}\\

{\bf Proof:} We fix the following curve in $\sneun$ joining the two elements $\pm 1 \in \sneun$:
$$ \gamma (t) = \cos (2t) + \sin (2t) v_1 \cdot v_2 \quad , \quad 0 \le t \le \frac{\pi}{2},$$

where $v_1 , \ldots , v_9$ is an orthonormal basis in $\rneun$. Using the explicit formulas for the spin representation $\kappa_9$ (see \cite{xx}) we compute the matrix of the endomorphism $\kappa_9 (\gamma (t)): \dneun \to \dneun$. $\kneun (\gamma(t))$ lifts into $Spin(16)$ via the formula
$$ \widetilde{\kneun (\gamma (t))}=( \cos (t) + \sin(t) e_1 \cdot e_2) \cdot \ldots \cdot (\cos(t) + \sin(t) e_{15} \cdot e_{16})$$

and at $t= \frac{\pi}{2}$ we obtain
$$ \widetilde{\kneun \Big(\gamma\Big( \frac{\pi}{2}\Big) \Big)} = e_1 \cdot e_2 \cdot \ldots \cdot e_{15} \cdot e_{16} . $$
\mbox{} \hfill \rule{3mm}{3mm}\\

Let us characterize $\sneun$ as a subgroup of $SO(16)$. For this purpose we consider the complex spin representation $\Delta^{\Bbb C}_9$ and recall that there exists a real structure $\alpha : \Delta^{\Bbb C}_9 \to \Delta^{\Bbb C}_9$ that anti-commutes with the Clifford multiplication of vectors by spinors (see \cite{xx}):
$$ \alpha (v \cdot \psi) = - v \alpha (\psi) \quad , \quad v \in {\Bbb R}^9, \psi \in \Delta^{\Bbb C}_9 . $$

The space $\dneun$ of real spinors is the fixed point set of $\alpha$:
$$ \dneun =\{ \psi \in \Delta^{\Bbb C}_9: \quad \alpha (\psi)= \psi \} . $$

We introduce a new multiplication of vectors by spinors via the formula
$$ v \star \psi := i(v \cdot \psi) . $$

 Since $\alpha$ is a real structure, the $\star$-multiplication is compatible with $\alpha$
$$ \alpha (v \star \psi) = v \star \alpha (\psi) . $$

In particular, the real spinor space $\dneun$ is invariant under the $\star$-multiplication by vectors. In case we understand the vectors of $\rneun$ as operators on $\dneun$ acting by the $\star$-multiplication, we will denote these vectors by $I, J \ldots \in \rneun$. Then we have\\

$(\star)$ \hfill $\displaystyle  I \star J + J \star I = 2 \langle I,J \rangle \quad , \quad I,J, \in \rneun . $ \hfill \mbox{}\\

Any vector $I \in \rneun$ defines a symmetric endomorphism $I: \dneun \to \dneun$ and, consequently, $\rneun$ is a subspace of the algebra $S_0^2 (\dneun)$ of all symmetric endomorphisms. For convenience, we will often omit the $\star$ between the vectors $I, J, \ldots \in {\Bbb R}^9$. \\

{\bf Proposition 2:} {\it  The group $\sneun$ consists of all products $(-1)^k I_1  I_2  \ldots I_{2k} \in SO(16)$ where $I_1, I_2, \ldots , I_{2k}$ are vectors of length one in ${\Bbb R}^9$. Moreover, the subgroup of $SO(\dneun)$ preserving under conjugation the space $\rneun \subset S_0^2 (\dneun)$ coincides with the group $\sneun$, i.e., 
$$ \sneun = \{ g \in SO(\dneun): \, \, g \, \rneun g^{-1} = \rneun \} . $$ }

{\bf Proof:} Consider the subgroup $ H= \{ g \in SO(\dneun): \, \, g \, \rneun g^{-1} = \rneun \} . $ Then we have 
$$ \sneun \subset H \subset SO(\dneun) . $$

On the other hand, $\sneun$ is a maximal nontrivial compact subgroup of $SO(16)$ (see \cite{x}). Consequently, we conclude $\sneun =H$. \mbox{} \hfill \rule{3mm}{3mm}\\

Let us consider  a 16-dimensional oriented Riemannian manifold $M^{16}$ and denote by ${\cal F}(M^{16})$ its frame bundle with structure group $SO(16)$.\\

{\bf Definition:} {\it A $\sneun$-structure is a reduction ${\cal R} \subset {\cal F}(M^{16})$ of the $SO(16)$-bundle ${\cal F}(M^{16})$ via the homomorphism $\kneun : \sneun \to SO(16)$}.\\

A $\sneun$-structure defines certain other geometric structures. In particular, it induces a spin structure on $M^{16}$ as well as a 9-dimensional real, oriented Euclidean vector bundle $V^9$ with spinor structure:
$$ V^9 := {\cal R} \times_{\sneun} \rneun . $$

The tangent bundle $T (M^{16})$ is isomorphic to the bundle $\dneun (V^9)$ of real spinors of the vector bundle $V^9$ and, therefore, we obtain a $\star$-multiplication of elements of $V^9$ by vectors in $T(M^{16})$. Conversely, a spin structure of $M^{16}$, a real vector bundle $V^9$ and a $\star$-multiplication define a $\sneun$-structure on $M^{16}$ (see Proposition 2). Locally a $\sneun$-structure is a collection of 9 sym\-metric \, involu\-tions $I_{\alpha} \, \, (1 \le \alpha \le 9)$ acting on the tangent bundle such that the following relations hold:
$$ I_{\alpha}^2 = \mbox{Id} \quad , \quad I^*_{\alpha} = I_{\alpha} \quad , \quad I_{\alpha} I_{\beta} = - I_{\beta} I_{\alpha} \quad (\alpha \not= \beta). $$

A $\sneun$-structure in dimension $n=16$ is the analogue of a quaternionic structure on Riemannian manifolds of dimension $n=4k$. The symmetric involutions $I_1 , \ldots , I_9$ define 2-forms $\Omega_{\alpha, \beta}$ on $M^{16}$ locally by the formula
$$ \Omega_{\alpha \beta} (X,Y) = g (X, I_{\alpha} I_{\beta} (Y)) - \langle I_{\alpha}, I_{\beta} \rangle g (X,Y) , \quad  X,Y \in T(M^{16}) . $$

The matrix $\Omega = (\Omega_{\alpha \beta})$ is an antisymmetric $(9 \times 9)$-matrix of 2-forms. Using the antisymmetric involutions $I_{\alpha}  I_{\beta}  I_{\gamma} \, \, (\alpha < \beta < \gamma)$ we can define in a similar way 2-forms $\Sigma_{\alpha \beta \gamma}$. Then
$$ \Sigma_{\alpha \beta \gamma} = - \Sigma_{\beta \alpha \gamma} = \Sigma_{\beta \gamma \alpha} $$

holds and the 2-forms $\{ \Omega_{\alpha \beta} , \Sigma_{\alpha \beta \gamma} \}$ are linearly independent and a local frame in the bundle $\Lambda^2 (M^{16})$.\\

The $\sneun$-representation $\Lambda^8 (\dneun) = \Lambda^8 ({\Bbb R}^{16})$ contains one and only one 8-form $\Omega_0^8$ which is invariant under the $\sneun$-action. This form defines the unique parallel form on the Cayley plane $F_4 /\sneun$. Since the signature of the Cayley plane is positive, $\Omega_0^8$ must be self-dual, $ \star \Omega_0^8 = \Omega_0^8$. It induces a {\it canonical 8-form} $\Omega^8$ defined on a 16-dimensional manifold $M^{16}$ with fixed $\sneun$-structure (see \cite{x}). 


\section{Large subgroups of $\sneun$}

\newcommand{\sacht}{Spin(8)}
\newcommand{\lneun}{\Lambda_9}
\newcommand{\lacht}{\Lambda_8}
\newcommand{\dam}{\Delta^-_8}

The group $\sacht$ admits an outer automorphism of order three (the principle of triality, see \cite{4x}, \cite{5x}). We use this automorphism to construct certain subgroups of $\sneun$ which are all pairwise not conjugate. In general, for any subgroup $H \subset \sneun$ we will denote by
\begin{itemize}
\item $\lneun (H)$ the representation of $H$ in the vector space $\rneun$;
\item $\dneun (H)$ the representation of $H$ in the spinor space $\dneun = {\Bbb R}^{16}$.
\end{itemize}

Therefore, we assign to any subgroup $H$ of $\sneun$  a pair $(\lneun (H), \dneun (H))$ of $H$-representations.\\

{\bf Example 1:} In case of the standard inclusion $\sacht \subset \sneun$ we have
$$ \lneun (\sacht)=\lacht \oplus {\Bbb R}^1 \quad , \quad \dneun (\sacht) = \Delta^+_8 \oplus \Delta^-_8  , $$

where $\lacht$ is the standard representation of $\sacht$ in ${\Bbb R}^8$ and $\Delta^{\pm}_8$ are the real spin representations of $\sacht$.\\

{\bf Example 2:} The kernel of the 8-dimensional real spin representation $\kappa^+_8 : \sacht \to SO(\Delta^+_8)$ is isomorphic to ${\Bbb Z}_2$:
$$ \ker (\kappa^+_8) =\{ 1, e_1 \cdot \ldots \cdot e_8 \} . $$

Consider the diagram
$$ \begin{diagram}
\node{\sacht} \arrow{s} \node{\mbox{} \hspace{1.4cm} \sacht \hspace{0.8cm} \subset } \arrow{s}   \node{\sneun} \arrow{s}\\
\node{\sacht / \ker (\kappa^+_8)} \arrow{e} \node{\mbox{} \, SO(\Delta^+_8) }  \node{SO(9)}\\
\end{diagram}$$
\vspace{-2.5cm}
$$ \mbox{} \hspace{4.2cm} \subset
$$

\vspace{0.5cm}

and lift the homomorphism
$$ \sacht \to \sacht / \ker (\kappa^+_8) \to SO(\Delta^+_8) $$

into the universal covering of $SO(\Delta^+_8)$. Then we obtain a subgroup $Spin^+(8) \subset \sneun$ isomorphic to $\sacht$. In this case we have

$$\lneun (Spin^+(8))= \Delta^+_8 \oplus {\Bbb R}^1 \quad , \quad \dneun (Spin^+(8))= \Delta_8^- \oplus \Delta^-_8 . $$

Indeed, the representation $\dneun$ splits under the action of the group $Spin^+(8)$ into two 8-dimensional irreducible representations. Moreover, the element $e_1 \cdot \ldots \cdot e_8 \in \sacht$ corresponds to the element $(-1) \in Spin^+(8) \subset \sneun$ and, consequently, $e_1 \cdot \ldots \cdot e_8$ acts on $\dneun$ by multiplication by $(-1)$. Therefore, we conclude that
$$ \dneun (Spin^+(8))= \Delta_8^- \oplus \Delta^-_8 . $$

{\bf Example 3:} The kernel of the 8-dimensional real spin representation $\kappa_8^- : \sacht \to SO(\Delta^-_8)$ is isomorphic to ${\Bbb Z}_2$
$$ \ker (\kappa^-_8) =\{ 1, - e_1 \cdot \ldots \cdot e_8 \} $$

and a construction similar to example 2 defines a subgroup $Spin^-(8) \subset \sneun$ isomorphic to $\sacht$ such that
$$ \lneun (Spin^-(8))= \Delta_8^- \oplus {\Bbb R}^1 \quad , \quad \dneun (Spin^-(8))= \Delta_8^+ \oplus \Delta_8^+ . $$

{\bf Example 4:}  The group $Spin(7)$ has a 7-dimensional real irreducible representation in ${\Bbb R}^7$ and an 8-dimensional real irreducible and faithful representation in the space $\Delta_7$ of real spinors. Lifting these two representations into $\sneun$ we obtain two subgroups $Spin(7)$ and $Spin_{\Delta} (7)$ of $\sneun$ such that
$$ \lneun (Spin(7))= \Lambda_7 \oplus {\Bbb R}^1 \oplus {\Bbb R}^1 \quad , \quad \dneun (Spin(7))= \Delta_7 \oplus \Delta_7 , $$
$$ \Lambda_9 (Spin_{\Delta}(7))= \Delta_7 \oplus {\Bbb R}^1 \quad , \quad \Delta_9 (Spin_{\Delta}(7))= \Delta_7 \oplus \Lambda_7 \oplus {\Bbb R}^1 . $$

The subgroup $Spin(7)$ is, in fact, already contained in $\sacht$. Consequently, we can once again apply the automorphism of the triality principle and obtain a total of 3 subgroups of $\sneun$, which we will denote by $Spin^+(7), Spin^-(7), Spin_{\Delta}(7)$.\\

{\bf Example 5:} We intersect the subgroup $\sneun \subset SO(16)$ with the subgroup $U(8) \subset SO(16)$. It turns out that $\sneun \cap U(8)$ is isomorphic to the group $Spin(2) \times_{{\Bbb Z}_2} Spin(7)$ and
$$ \dneun (\sneun \cap U(8))= {\Bbb R}^2 \otimes \Delta_7 \quad , \quad \Lambda_9 (\sneun \cap U(8))= {\Bbb R}^2 \oplus \Lambda_7 . $$

Table 1 summarizes  the decomposition of the representations $\Lambda_9$ and $\dneun$ for all these subgroups.\\

\begin{center}
\begin{tabular}{|c|c|c|}
\hline &&\\
$H$ & $\lneun (H)$ & $\dneun (H)$\\ &&\\ \hline
&&\\
$\sacht$ & $\Lambda_8 \oplus {\Bbb R}^1$ & $\Delta_8^+ \oplus \Delta_8^-$\\
&&\\ \hline &&\\
$Spin^+(8)$ & $\Delta_8^+ \oplus {\Bbb R}^1$ & $\Delta_8^- \oplus \Delta_8^-$\\
&&\\ \hline &&\\
$Spin^-(8)$ & $\Delta_8^- \oplus {\Bbb R}^1$ & $\Delta_8^+ \oplus \Delta_8^+$\\
&&\\ \hline &&\\
$\sneun \cap U(8)$ & ${\Bbb R}^2 \oplus \Lambda_7$ & ${\Bbb R}^2 \otimes \Delta_7$\\
&&\\ \hline &&\\
$Spin(7)$ & $\Lambda_7 \oplus {\Bbb R}^1 \oplus {\Bbb R}^1$ & $\Delta_7 \oplus \Delta_7$\\
&&\\ \hline &&\\
$Spin_{\Delta}(7)$ & $\Delta_7  \oplus {\Bbb R}^1$ & $\Lambda_7 \oplus \Delta_7 \oplus {\Bbb R}^1$\\
&&\\ \hline &&\\
$Spin^+(7)$ & $\Delta_7 \oplus {\Bbb R}^1$ & $\Delta_7 \oplus \Delta_7$\\
&&\\ \hline &&\\
$Spin^-(7)$ & $\Delta_7 \oplus {\Bbb R}^1$ & $\Delta_7 \oplus \Delta_7$\\
&&\\ \hline &&\\
$G_2$ & $\Lambda_7 \oplus {\Bbb R}^1 \oplus {\Bbb R}^1$ & $\Lambda_7 \oplus \Lambda_7 \oplus {\Bbb R}^1 \oplus {\Bbb R}^1$.\\ &&\\ \hline
\end{tabular}

\vspace{0.5cm}

{\small Table 1: Large subgroups of $\sneun$ and the branching of the standard and \\
\mbox{} \hspace{-4.3cm} spin representation of $\sneun$.}\\
\end{center}

Since the subgroup $Spin_{\Delta}(7) \subset \sneun$ is the isotropy group of a spinor in $\dneun$, we have\\

{\bf Proposition 1:} {\it Let $M^{16}$ be an oriented Riemannian manifold admitting a non-vanishing vector field. Then $M^{16}$ admits a $\sneun$-structure if and only if $M^{16}$ admits a $Spin_{\Delta} (7)$-structure. In case $M^{16}$ admits two independent vector fields, the existence of a $\sneun$-structure is equivalent to the existence of a $G_2$-structure.}\\


\section{Examples of $\sneun$-structures}

A $\sneun$-structure on a 16-dimensional Riemannian manifold occurs in case the frame bundle admits a reduction to a subgroup of $\sneun \subset SO(16)$. We describe two cases of a geometric situation of this type. First, we will construct examples of homogeneous spaces with invariant $\sneun$-structures. Consider a Lie group $G$, a subgroup $H \subset G$ and suppose that the homogeneous space $G/H$ is reductive. We decompose the Lie algebra $\g$ of $G$ into
$$ \g = \h \oplus \n , $$

where the subspace $\n$ is $Ad(H)$-invariant. If $H$ is a subgroup of $\sneun$ and the $H$-representations $\dneun (H)$ and $Ad:H \to SO (\n)$ are equivalent, then the space $G/H$ admits a homogeneous  $\sneun$-structure. Indeed, the frame bundle ${\cal F}(G/H)$ given by
$$ {\cal F}(G/H) = G \times_{Ad} SO (\n) $$

admits a $H$-reduction and the subgroup $H$ is contained in $\sneun$. This general remark yields the following examples of 16-dimensional manifolds with homogeneous $\sneun$-structures.\\

{\bf Example 1:} Take $G=F_4$ or $F_4^*$ and  $H= \sneun$. Then the  symmetric spaces $F_4 / \sneun$ and $F_4^*/ \sneun$ admit homogeneous $\sneun$-structures.\\

{\bf Example 2:} Consider the subgroup $H= Spin_{\Delta}(7) \subset \sneun$. Then we have
$$ \dneun (Spin_{\Delta}(7))= \Delta_7 \oplus \Lambda_7 \oplus {\Bbb R}^1 . $$

On the other hand, $\sneun$ acts transitively on the sphere $S(\Delta_9) = S^{15}$ with isotropy group $Spin_{\Delta}(7)$. The isotropy representation of the homogeneous space $\sneun / Spin_{\Delta}(7)$ is isomorphic to $\Delta_7 \oplus \Lambda_7$. Consequently, the pair of groups
$$ G=S^1 \times \sneun \quad , \quad H= Spin_{\Delta} (7) $$
defines a homogeneous $Spin_{\Delta}(7) \subset \sneun$-structure on the manifold $ G/H = S^1 \times S^{15} . $\\

{\bf Example 3:} Consider the pair
$$ G=S^1 \times S^1 \times SO(8) \quad , \quad H=G_2 , $$

where $G_2$ denotes the exceptional Lie group embedded into $SO(7) \subset SO(8)$. The isotropy representation of this homogeneous space is isomorphic to ${\Bbb R}^1 \oplus {\Bbb R}^1 \oplus \Lambda_7 \oplus \Lambda_7$ and coincides with the $G_2$-representation $\Delta_9(G_2)$. Consequently, the homogeneous space
$$ G/H =S^1 \times S^1 \times (SO(8)/G_2)$$

admits a homogeneous $G_2 \subset \sneun$-structure.\\

{\bf Example 4:} The group $H=SU(3)$ is the isotropy group of a pair of real spinors in $\Delta_7$ and, henceforth, a subgroup of $Spin(7) \subset \sneun$. Consider the group $G= S^1 \times S^1 \times S^1 \times Spin(7)$. The isotropy representation of the homogeneous space $G/H$ is isomorphic to ${\Bbb R}^1 \oplus {\Bbb R}^1 \oplus {\Bbb R}^1 \oplus {\Bbb R}^1 \oplus {\Bbb R}^6 \oplus {\Bbb R}^6$ and coincides with the $H$-representation $\dneun (SU(3))$. Consequently, the homogeneous space
$$ G/H = S^1 \times S^1 \times S^1 \times (Spin(7)/ SU(3)) $$

admits a homogeneous $SU(3) \subset Spin(9)$-structure.\\

{\bf Example 5:} The pair $ G= SU(5) \quad , \quad H= SU(3)$ defines a homogeneous $SU(3) \subset \sneun$-structure on
$$ G/H = SU(5)/ SU(3) . $$

The second situation in which a $\sneun$-structure occurs in a natural way is the case where the tangent bundle of $M^{16}$ splits  in a suitable way.\\

{\bf Proposition 1:} {\it Let $M^{16}$ be a 16-dimensional oriented Riemannian manifold and suppose that there exist a 4-dimensional complex vector bundle $E^4$ as well as a complex line bundle $L$ such that
\begin{enumerate}
\item $c_1 (E^4)=0$ \, \, in \, \, $H^2 (M^{16}; {\Bbb Z})$;
\item the tangent bundle $T(M^{16})$ is, as a real vector bundle, isomorphic to the Whitney sum $L \otimes (E^4 \oplus E^4)$.
\end{enumerate}

Then $M^{16}$ admits a $\sneun \cap U(8) \subset \sneun$-structure.}\\

{\bf Proof:} Since $SU(4)$ is isomorphic to $Spin(6)$, the tangent bundle of $M^{16}$ admits a $Spin(2) \times_{{\Bbb Z}_2} Spin(6) \subset \sneun \cap U(8)$-reduction. The representation $\dneun (Spin(6))= \Delta_9 (SU(4))$ is then isomorphic to the standard representation of $SU(4)$ in ${\Bbb C}^4 \oplus {\Bbb C}^4$. \hfill \rule{3mm}{3mm}\\

{\bf Proposition 2:} {\it Let $M^{16}$ be a 16-dimensional oriented Riemannian spin manifold. Suppose that there exists a 8-dimensional real vector bundle $W^8$ with $Spin(7)$-structure such that the tangent bundle $T(M^{16})$ is isomorphic to $W^8 \oplus W^8$. Then $M^{16}$ admits a $Spin(7) \subset \sneun$-structure.}\\

{\bf Proof:} Consider the subgroup $H= Spin(7)$. Then we know already that $ \dneun (H)= \Delta_7 \oplus \Delta_7 $ holds. \hfill \rule{3mm}{3mm}\\

Let $N^k$ be an arbitrary manifold and consider the projection $\pi :T(N^k) \to N^k$ of its tangent bundle. The bundle $T(T(N^k))$ is isomorphic to the sum of the induced bundles $\pi^* (T(N^k)) \oplus \pi^* (T(N^k))$. This isomorphism is not a canonical one, but depends on a fixed linear connection on the manifold $N^k$. Therefore, we obtain the following\\

{\bf Corollary 1:} {\it Let $N^8$ be an oriented, 8-dimensional Riemannian manifold with a $Spin(7)$-structure. Then the tangent bundle $M^{16}=T(N^8)$ admits a $Spin(7) \subset \sneun$-structure.}\\

{\bf Corollary 2:} {\it Let $N^8$ be an 8-dimensional Hermitian manifold with first Chern class divisible by 4. Then the tangent bundle $M^{16} =T(N^8)$ admits a $\sneun \cap U(8) \subset  \sneun$-structure.}\\

\section{Topological conditions}

A 16-dimensional compact manifold with $\sneun$-structure should satisfy certain topological conditions. Some of them have already been studied in the paper \cite{5x} and we will first summarize these results.\\

\newcommand{\vn}{V^9}
\newcommand{\om}{\omega}
\newcommand{\ms}{M^{16}}

{\bf Theorem 1 (see \cite{5x}):} {\it Let $M^{16}$ be a compact manifold admitting a $\sneun$-structure and denote by $V^9$ the associated 9-dimensional bundle. 

\begin{enumerate}
\item The following Stiefel-Whitney classes of $M^{16}$ vanish
$$ \om_1 = \om_2 = \om_3 = \om_4 = \om_5 = \om_6 = \om_7 = 0$$
$$ \om_9 =\om_{10} = \om_{11} = 	\om_{13} =0 . $$
\item The Stiefel-Whitney classes of $M^{16}$ are related to the corresponding classes of $\vn$ by the formulas
\begin{eqnarray*}
\om_8 (\ms) &=& \om^2_4 (\vn) + \om_8 (\vn)\\
\om_{12} (\ms) &=& \om_6^2 (\vn) + \om_4 (\vn) \, \om_8 (\vn)\\
\om_{14} (\ms) &=& \om_7^2 (\vn) + \om_6 (\vn) \, \om_8 (\vn)\\
\om_{15} (\ms) &=& \om_7 (\vn) \, \om_8 (\vn) . 
\end{eqnarray*}
\item In case $H^* (\ms; {\Bbb Z})$ is 2-torsion free, the Pontrjagin class $p_1 (\ms)$ is divisible by 4 and the Pontrjagin classes $p_2 (\ms) \,  ,  \,  p_3 (\ms)$ are divisible by 2.
\end{enumerate}}

Since the first seven Stiefel-Whitney classes of $M^{16}$ vanish, the Wu class of $M^{16}$ reduces to the element $\omega_{8} {(M^{16})} \in H^8 (M^{16}; {\Bbb Z}_2)$ (see \cite{huse}). Consequently, the Stiefel-Whitney class $\omega_8 (M^{16})$ is characterized by the condition
$$ y^8 \cup y^8 = y^8 \cup \omega_8 (M^{16}) \quad \mbox{for any $y^8 \in H^8 (M^{16}; {\Bbb Z}_2)$}  . $$

{\bf Corollary 1:} {\it Let $M^{16}$ be a compact manifold admitting a $\sneun$-structure. Then the quadratic form over ${\Bbb Z}$
$$ H^8 (M^{16}; {\Bbb Z}) / \mbox{Tor} \times H^8 (M^{16}; {\Bbb Z}) / \mbox{Tor} \to H^{16} (M^{16}; {\Bbb Z}) $$

is an even ${\Bbb Z}$-form if and only if $\omega_8 (M^{16})=0$.}\\

The aim of this part of the paper is to compute the Pontrjagin classes of $\ms$ explicitly in terms of the corresponding classes of the vector bundle $\vn$. Some new integral conditions are consequences of these formulas.\\

{\bf Theorem 2:} {\it Let $\ms$ be a compact manifold admitting a $\sneun$-structure. Then the Pontrjagin classes of $\ms$ and of the bundle $\vn$ are related by the following formulas:
\begin{enumerate}
\item $p_1 (\ms)= 2 \, p_1 (\vn)$.
\item $p_2 (\ms)= \frac{7}{4}  \, p_1^2 (\vn) - p_2 (\vn)$.
\item $p_3 (\ms)= \frac{1}{8} \Big( 7 \, p_1^3 (\vn) - 12 \, p_1 (\vn) \, p_2 (\vn) + 16 \, p_3 (\vn)\Big).$
\item $p_4 (\ms) = \frac{1}{128} \Big(35 \, p_1^4 (\vn) - 120 \, p_1^2 (\vn) \, p_2 (\vn) + 400 \, p_1 (\vn) \, p_3 (\vn) - 1664 \, p_4 (\vn)\Big).$
\end{enumerate}

The Euler class $e(M^{16})$ and the the fourth $L$-polynomial of $\ms$ can be expressed by the Pontrjagin classes of the bundle $V^9$:
\begin{itemize}
\item[5.] $\displaystyle e(M^{16})= \frac{1}{256} \, p^4_1 (V^9) - \frac{1}{32} \, p_1^2 (V^9) p_2 (V^9) + \frac{1}{16} \, p_2^2 (V^9) - \frac{1}{4} \, p_4 (V^9) . $
\item[6.] $\displaystyle L_4  (\ms) = \frac{1}{1814400} \Big( 3551 \, p_1^4 (\vn) - 21208 \, p_1^2 (\vn) \, p_2 (\vn) $

\smallskip
\mbox{} \hspace{1.8cm} $\displaystyle  + 116048 \, p_1 (\vn) \, p_3 (\vn) - 128 ( 19 \, p_2^2 (\vn) + 4953 \, p_4 (\vn))\Big) . $\\
\end{itemize}}

The bundle $\vn$ is a real vector bundle with spin structure. The general formula
$$ p_1 (\vn) \equiv \om_2^2 (\vn) =0 \quad \mbox{mod} \, 2 $$

yields that $p_1 (\vn)$ is divisible by 2. In particular, we obtain the\\

{\bf Corollary 2:} {\it Let $\ms$ be a compact manifold admitting a $\sneun$-structure. Then,
\begin{enumerate}
\item $p_1 (\ms)$ is divisible by 4. Denote by $\mbox{\normalsize x}$ the cohomology class ${p_1 (\ms)}/{4}$.
\item $\frac{1}{2} (p_3 (\ms) - 3 \, \mbox{\normalsize x} \, p_2 (\ms))$ is an integral cohomology class.
\item $\frac{175}{8} \, \mbox{\normalsize x}^4 - \frac{45}{8} \, \mbox{\normalsize x}^2 \, p_2 (\ms) + \frac{25}{8} \, \mbox{\normalsize x} \, p_3 (\ms) - p_4 (\ms)$ is divisible by 13.
\end{enumerate}}

\medskip

{\bf Proof of Theorem 2:} Consider the 9-dimensional spin representation
$$ \kappa_9 : \sneun \to SO(\dneun) $$

and fix maximal tori $T^4, T^8$ in $SO(9)$ and $SO(\dneun)$. We denote by \linebreak $\Theta_i : T^4 \to S^1 \, \, (1 \le i \le 4)$ and by $\mu_{\alpha} : T^8 \to S^1 \, \, (1 \le \alpha \le 8)$ the coordinates of the maximal tori. Then, the representation $\kappa_9$ has the following weights:\\

\begin{tabular}{ll}
$\displaystyle \mu_1 = \frac{1}{2} (\Theta_1 + \Theta_2 + \Theta_3 + \Theta_4) \quad , $&$\displaystyle \mu_2 = \frac{1}{2} (\Theta_1 + \Theta_2 + \Theta_3 - \Theta_4)$\\
&\\
$\displaystyle \mu_3 = \frac{1}{2} (\Theta_1 + \Theta_2 - \Theta_3 + \Theta_4) \quad , $&$\displaystyle  \mu_4 = \frac{1}{2} (\Theta_1 - \Theta_2 + \Theta_3 + \Theta_4)$\\
&\\
$\displaystyle \mu_5 = \frac{1}{2} (- \Theta_1 + \Theta_2 + \Theta_3 + \Theta_4) \, \, \, , $&$\displaystyle  \mu_6 = \frac{1}{2} (\Theta_1 + \Theta_2 - \Theta_3 - \Theta_4)$\\
&\\
$\displaystyle \mu_7 = \frac{1}{2} (\Theta_1 - \Theta_2 + \Theta_3 - \Theta_4) \quad , $&$\displaystyle  \mu_8 = \frac{1}{2} (- \Theta_1 + \Theta_2 + \Theta_3 - \Theta_4)$\\
\end{tabular}

\vspace{0.5cm}

The first Pontrjagin class $p_1 (\ms)$ is given by
$$ p_1 (\ms) = \sum\limits^8_{\alpha=1} \mu_{\alpha}^2 = 2 \sum\limits^4_{i=1} \Theta_i^2 = 2 p_1 (\vn) . $$

This calculation proves the first formula of Theorem 2. For the second Pontrjagin class we obtain
$$p_2 (\ms) = \sum\limits_{1 \le \alpha < \beta \le 8} \mu_{\alpha}^2 \mu_{\beta}^2 =  \frac{7}{4} \Big( \sum\limits_{i=1}^4 \Theta_i^2 \Big)^2 - \sum\limits_{1 \le i < j \le 4} \Theta_i^2 \Theta_j^2 = \frac{7}{4} p_1^2 (\vn) - p_2 (\vn) . 
$$

The formulas 3.) - 5.) can be computed in a similar way, however the calculations of $p_3, \, p_4$ are much more lengthly. 
The last formula  is a consequence of the first four formulas and the formula for the $L_4$-polynomial
$$ L_4 = \frac{1}{3^4 \, 5^2 \, 7} (381 \, p_4 - 71 \, p_3 \, p_1 - 19 \, p_2^2 + 22 \, p_2 \, p_1^2 - 3 \, p_1^4) , $$
(see \cite{7x}, page 25).  \hfill \rule{3mm}{3mm}\\

{\bf Example:} Let $\ms$ be a smooth complete intersection of three quadrics in ${\Bbb P}^{11} ({\Bbb C})$. Then the diffeomorphism type is unique and $\ms$ is a simply-connected manifold (see \cite{8x}). Denote by $\mbox{\normalsize x} \in H^2 ({\Bbb P}^{11} ({\Bbb C}); {\Bbb Z})$  the generator of the second cohomology group. The Stiefel-Whitney, Chern and Pontrjagin classes of $\ms$ are well-known (see \cite{7x}, page 159):\\

\mbox{} \hspace{2cm} $\omega_i=0 \quad \mbox{for $ i \not= 8, \, \, \, \omega_8 = x^4$}$\medskip

\begin{center}
\begin{tabular}{lllllll}
$\displaystyle c_1=6 \, \mbox{\normalsize x}$ & , & $c_2 =18 \, \mbox{\normalsize x}^2$ &  , &  $c_3 = 32 \, \mbox{\normalsize x}^3$ &  , &  $c_4 =39 \, \mbox{\normalsize x}^4$  ,  \\
&&&&&&\\
$c_5=30 \, \mbox{\normalsize x}^5$ &  , &  $c_6 =20 \, \mbox{\normalsize x}^6$ &  , &  $c_7 = 0$ &  , &  $c_8 = 15 \, \mbox{\normalsize x}^8 $ , \\
&&&&&&\\
$p_1=0$ &  , &  $p_2 =18 \, \mbox{\normalsize x}^4$ &  , &  $p_3 = 60\,  \mbox{\normalsize x}^6$ &  , &  $p_4 = 351 \, \mbox{\normalsize x}^8$ . \\ 
\end{tabular}
\end{center}

\vspace{0.3cm}

The inclusion $M^{16} \to {\Bbb P}^{11} ({\Bbb C})$ induces an isomorphism $H_i (\ms; {\Bbb Z}) \to H_i ({\Bbb P}^{11} ({\Bbb C}); {\Bbb Z})$ for $i <8$ (see \cite{grady}). In particular, $\mbox{x} \in H^2 (\ms; {\Bbb Z})$ is a generator of the second homology group of $\ms$ and $\mbox{x}^4 \in H^8 (\ms; {\Bbb Z})$ is not divisible.   A calculation of the Euler class as well as the signature yields the equality
$$ \frac{{\cal X}(\ms)}{\sigma (\ms)} = \frac{5}{3} . $$

In particular, $\ms$ satisfies all the necessary conditions for the Pontrjagin classes formulated in Corollary 1 and Corollary 2. Therefore, complete intersections of three quadrics in ${\Bbb P}^{11} ({\Bbb C})$ are candidates of compact 16-dimensional manifolds admitting $\sneun$-structures. However, a $\sneun$-structure compatible with the complex structure on $\ms$ cannot exist. Indeed, since $\sneun \cap U(8)$ is isomorphic to $Spin(2) \times_{{\Bbb Z}_2} Spin(7)$ and $Spin(7)$ is contained in the subgroup $SU(8) \subset SO(16)$, such a structure would define a product decomposition of the complex tangent bundle $T(\ms)$ into $T(M^{16}) = L \otimes F^8$, where $L$ is an complex line bundle and $F^8$ is an 8-dimensional complex vector bundle with vanishing first Chern class, $c_1 (F^8)=0$. Then the first Chern class of $\ms$ is divisible by 8, $c_1 (\ms)=8 \, c_1 (L)$, a contradiction. \\

These formulas become much simpler in case  $p_1 (\ms)=0$ vanishes and the scalar curvature $R>0$ is positive.\\

{\bf Corollary 3:} {\it Let $\ms$ be a compact 16-dimensional Riemannian manifold with $\sneun$-structure and suppose that the scalar curvature $R>0$ is positive. Moreover, let $p_1 (\ms)=0$ be trivial. Then,
\begin{enumerate}
\item the signature $\sigma (\ms)$ is given by the formula
$$ \sigma (\ms) = - \frac{1}{3} \int\limits_{\ms} p_4 (\vn) = \frac{1}{39} \int\limits_{\ms} p_4 (\ms) .$$
\item in $H^*(\ms; {\Bbb Q})$ the formulas $p_2 (\ms) = - \, p_2 (\vn), \, \, p_3 (\ms) = 2 \, p_3 (\vn)$ and \, \, $p_4 (\ms) = - 13 \, p_4 (\vn)$ hold.
\item the equality $13 \, p_2^2 (\ms) = 12 \, p_4 (\ms)$ holds.
\item the Euler class is given by $e(\ms)= \frac{1}{16} (p_2^2(\vn) - 4 \, p_4 (\vn))= \frac{1}{13} \, p_4 (\ms)$. In particular, we obtain $e(\ms) = 3 \, \sigma (\ms)$.
\end{enumerate}}

{\bf Proof of Corollary 3:} $\ms$ does not admit harmonic spinors, and, consequently, the $\hat{\cal A}$-genus vanishes. Since $p_1 (\ms)=0$, the $\hat{\cal A}$-genus is a combination of two characteristic classes only: 
$$ \hat{\cal A} = A \, p_2^2 + B \, p_4 . $$

$A$ and $B$ are universal coefficients. In the case of the Cayley plane $\ms = F_4 / \sneun$, we know that
$$ p_2^2 = 36 \quad , \quad p_4 = 39 ,$$

(see \cite{6x}). Consequently, the ratio $A/B$ equals $- \frac{39}{36}$ and, finally, we obtain the equation
$$ 39 \, p_2^2 (\ms) = 36 \, p_4 (\ms)  $$

for any manifold $\ms$ under consideration. The signature of $\ms$ is given by
\begin{eqnarray*}
\sigma (\ms) &=& \int\limits_{\ms} L_4 \, \, =  - \frac{1}{14175} \int\limits_{\ms} (19 \, p_2^2 (\vn) + 4953 \, p_4 (\vn))\\
&=&  - \frac{1}{14175} \int\limits_{\ms} (19 \, p_2^2 (\ms) + 4953 \, p_4 (\vn))\\
&=&  - \frac{1}{14175} \int\limits_{\ms} \Big(19 \cdot \frac{36}{39} \, p_4 (\ms) + 4953 \, p_4 (\vn)\Big)\\
&=&   - \frac{1}{14175} \int\limits_{\ms} \Big(19 \cdot \frac{36}{39} \cdot \Big( - \frac{1664}{128} \Big) + 4953 \Big) \, p_4 (\vn) = - \frac{1}{3} \int\limits_{\ms} \, p_4 (\vn) .
\end{eqnarray*}
\mbox{} \hfill \rule{3mm}{3mm}\\

\section{The homotopy type of the space $Spin (16)/ Spin (9)$}

A necessary and sufficient criterion for the existence of  a $\sneun$-structure does not seem to be known. The classifying space of the group $\sneun$ is a fibre bundle over the classifying space of $Spin (16)$ with fibre $X^{84} = Spin(16) / \sneun$. Therefore, the homotopy  type of $X^{84}$ yields obstructions for the existence of a $\sneun$-structure on a 16-dimensional real vector bundle. We compute some of the homotopy groups of this space. \\

\newcommand{\xav}{(X^{84})}

{\bf Theorem 1:} {\it
\begin{enumerate}
\item $\pi_1 \xav = \pi_2 \xav = \pi_4 \xav = \pi_5 \xav = \pi_6 \xav = \pi_7 \xav = \pi_{13} \xav =0.$
\item $\pi_3 \xav = {\Bbb Z}_4 \quad , \, \, \, \pi_{12} \xav = \pi_{14} \xav = {\Bbb Z}_2.$
\item There is an exact sequence 
$$  0 \to \pi_{10} \xav  \to {\Bbb Z}_2 \oplus {\Bbb Z}_2 \to {\Bbb Z}_2 
\to  \pi_9 \xav \to {\Bbb Z}_2 \oplus {\Bbb Z}_2 \to {\Bbb Z}_2 
\to \pi_8 \xav \to 0 . 
$$
\item There is a surjective homomorphism $\pi_{11} \xav \to {\Bbb Z}_8$.
\end{enumerate}}

Using this result as well as the condition for the first Pontrjagin class discussed before we immediately obtain the following\\

{\bf Corollary:} {\it Let $Y$ be an 8-dimensional $CW$-complex. A real, oriented 16-dimensional vector bundle over $Y$ admits a $\sneun$-structure if and only if its first Pontrjagin class is divisible by 4.}\\

{\bf Proof of the Theorem:} We apply the exact sequence of homotopy groups of the fibration
$$ \sneun \to Spin(16) \to X^{84} . $$

Since $\pi_{12} (Spin(16))= \pi_{13} (Spin(16))= \pi_{14}(Spin(16))=0$, we obtain
$$ \pi_{14} \xav = \pi_{13} (\sneun) = {\Bbb Z}_2 ,$$
$$ \pi_{13} \xav = \pi_{12} (\sneun)=0  , $$

(see \cite{9x}). Next we use the fact that $\pi_{12} (Spin(16)) = \pi_{10} (Spin(16))=0$. Then, we obtain the exact sequence
$$ 0 \to \pi_{12} \xav \to \pi_{11} (\sneun) \to \pi_{11} (Spin(16)) \to \pi_{11} \xav \to \pi_{10} (\sneun) \to 0 . $$

But $\pi_{11} (\sneun) = {\Bbb Z} \oplus {\Bbb Z}_2$, $\pi_{10} (\sneun)= {\Bbb Z}_8$ and $\pi_{11} (Spin(16)) = {\Bbb Z}$ \, (see \cite{9x}). Consequently, we obtain $\pi_{12} \xav = {\Bbb Z}_2$ as well as a surjective homomorphism $\pi_{11} \xav \to {\Bbb Z}_8$. The other statements of the theorem are easy consequences of the following two facts, which we are going to prove now:
\begin{itemize}
\item[a.)] The homomorphism induced by the inclusion $\kappa_9 : \sneun \to Spin(16)$
$$ (\kneun)_{\#} : \pi_3 (\sneun) = {\Bbb Z} \to {\Bbb Z} = \pi_3 (Spin(16))  $$

is the multiplication by 4; 
\item[b.)] The inclusion $\kneun$ induces an isomorphism
$$ (\kneun)_{\#} : \pi_7 (\sneun) \to \pi_7 (Spin(16)) . $$
\end{itemize}

Indeed, consider the subgroup $Spin(3) \subset \sneun$. The homogeneous space $\sneun / Spin(3)$ is the Stiefel manifold $V_3 ({\Bbb R}^9)$ and its homotopy groups are well-known:
$$ \pi_3 (V_3({\Bbb R}^9))= \pi_4 (V_3({\Bbb R}^9))=0 . $$

Therefore, the inclusion induces an isomorphism $\pi_3 (Spin(3)) \to \pi_3 (\sneun)$. Since $\dneun (Spin(7)) = 2 \Delta_7$, we can calculate $\Delta_9 (Spin(3))$ by restricting the 7-dimensional spin representation $\Delta_7$ to $Spin(3)$. But this restriction coincides with
$$ \Delta_{7} | Spin(3) = {\Bbb C}^2 \oplus {\Bbb C}^2 , $$

where $Spin(3) = SU(2) \subset SO(4)$ acts on ${\Bbb C}^2$ in the usual way. Finally, we obtain a.).\\

The proof of the property b.) is more sophisticated and uses some  results of \cite{3x}. Consider the subgroup
$$ Spin(7) \cong Spin_{\Delta} (7) \subset \sneun . $$

Since $\sneun$ acts transitively on $S^{15}$ with isotropy group $Spin_{\Delta} (7)$, the homomorphism
 $$ \pi_7 (Spin_{\Delta} (7)) \to \pi_7 (\sneun) $$

is an isomorphism. Because of $\dneun (Spin_{\Delta}(7))= \Delta_7 \oplus \Lambda_7 \oplus {\Bbb R}^1$ we should study the map
$$ Spin(7) \to SO(\Delta_7) \oplus SO(\Lambda_7) . $$

Let $i : Spin(7) \to Spin(8)$ be the inclusion and denote by $H: Spin(8) \to Spin(8)$ the triality automorphism. The 7-dimensional spin representation $\kappa_7$ is given by $\kappa_7 = H \circ i$. Denote by $p: Spin(8) \to S^7$ the projection and fix a generator $\alpha_7 \in \pi_7 (Spin(7))= {\Bbb Z}$. Moreover, we choose generators $e_1, e_2 \in \pi_7 (Spin(8)) = {\Bbb Z} \oplus {\Bbb Z}$ such that 
\begin{enumerate}
\item[1.)] $p_{\#} (e_1) = e $ is a generator of the group $\pi_7 (S^7)$;
\item[2.)] $i_{\#} (\alpha_7) = e_2.$ 
\end{enumerate}

There are elements $\alpha, \beta, \gamma \in \pi_7 (Spin(8))$ with the following properties (see \cite{3x}, page 152-155):
\begin{itemize}
\item[a.)] $\gamma = \alpha - \beta$;
\item[b.)] $\alpha = e_1 + xe_2  , \quad \beta = - e_1 + ye_2  , \quad \gamma = 2e_1 + (x-y) e_2$ \, where $x,y \in {\Bbb Z}$ are integers;
\item[c.)] The element $\gamma$ maps into $0 \in \pi_7 (\sneun)$ via the homomorphism $\pi_7 (Spin(8)) \to \pi_7 (Spin(9))$;
\item[d.)] The triality homomorphism acts on $\pi_7 (Spin(8))$ via the following formulas:
$$ H_{\#} (\alpha)= - \beta \quad , \quad H_* (\beta) = \gamma \quad , \quad H_{\#} (\gamma)= - \alpha . $$
\end{itemize}

Let 
$$ H_{\#} = \left( \begin{array}{cc} A&C\\B&D \end{array} \right) $$

be the matrix of $H_{\#} : \pi_7 (Spin(8)) \to \pi_7 (Spin(8))$ with respect to the basis $e_1, e_2$. Then, condition d.) is equivalent to the system of six equations:
$$ A +xC=1 \quad , \quad - A +yC=2 \quad , \quad 2A + (x-y)C= - 1$$
$$B+xD = - y \quad , \quad - B+yD =x-y \quad , \quad 2B+(x-y)D= - x.$$

We solve this system:
$$ C = - \frac{1+A+A^2}{B} \quad , \quad D=-1-A \quad , \quad x= \frac{AB-B}{1+A+A^2} \quad , \quad y =  - \frac{2B+AB}{1+A+A^2} $$

and, in particular, we obtain
$$ x-y =C^{-1} (1+2D) . $$

Consequently, we have $\gamma = 2e_1 + {C^{-1}} {(1+2D)} e_2$. On the other hand, $\gamma=0$ in $\pi_7(Spin(9))$, and, therefore, we conclude that $ 2e_1 = - C^{-1} (1+2D)  e_2 $ holds in $\pi_7 (Spin(9))$. This implies $ H_{\#} (e_2) = C e_1 + De_2 = - \frac{1}{2} e_2 $ in $\pi_7 (\sneun)$. This equation implies that the homomorphism
$$ (\kappa_7)_{\#} : \pi_7 (Spin(7)) \to \pi_7 (SO(\Delta_7)) \to \pi_7 (\sneun)$$

maps the generator $\alpha_7 \in \pi_7 (Spin(7))$ into $(- \alpha_9)$, where $\alpha_9 \in \pi_7 (\sneun)$ is the generator of $\pi_7 (\sneun)$. The inclusion $Spin(7) \to \sneun$ induces the map $\alpha_7 \to 2 \alpha_9$. Indeed, the homotopy group
$$ \pi_7 (\sneun / Spin(7)) = \pi_7 (V_{2,9})= {\Bbb Z}_2$$

is isomorphic to ${\Bbb Z}_2$. Finally, we conclude that the map
$$ Spin(7) \to SO(\Delta_7) \oplus SO(\Lambda_7) \subset SO(16)$$

induces an isomorphism on $\pi_7$. \hfill \rule{3mm}{3mm}\\

\section{The decomposition of $\sneun$-representations}

Let us once again recall the notation. An orthonormal basis in ${\Bbb R}^9$ is denoted by $I_1, \ldots , I_9$. Via the modified Clifford multiplication $\star$ (see Section 2) the vectors $I_{\alpha} \, \, (1 \le \alpha \le 9)$ are symmetric involutions acting on the space $\Delta_9 = {\Bbb R}^{16}$:
$$ I_{\alpha}^2 = \mbox{Id} \quad , \quad I_{\alpha}^* = I_{\alpha} \quad , \quad I_{\alpha}  I_{\beta} = - I_{\beta}  I_{\alpha} \quad (\alpha \not= \beta) . $$

The group $\sneun$ acts on the space $\Lambda^k (\rneun) \otimes \dneun$ of spinor valued $k$-forms in $\rneun$. The decomposition into irreducible components of this space is well-known  (see  \cite{10x}). Let us introduce the endomorphisms
$$ \Theta_k : \Lambda^k (\rneun) \otimes \dneun \to \Lambda^{k+1} (\rneun) \otimes \dneun , $$
$$ \Theta_k^* : \Lambda^k (\rneun) \otimes \dneun \to \Lambda^{k-1} (\rneun) \otimes \dneun , $$

defined by the formulas
$$ \Theta_k (\om^k \otimes \varphi)= \sum\limits^9_{\alpha=1} (I_{\alpha} \wedge \om^k) \otimes (I_{\alpha} \star \varphi) , $$
$$ \Theta_k^* (\om^k \otimes \varphi)=  - \sum\limits^9_{\alpha=1} (I_{\alpha} \, \mbox{}_{\_\!\!\_\!}{\scriptstyle{\mid}} \, \om^k) \otimes (I_{\alpha} \star \varphi)  $$

and denote by $P_{r}$ the kernel of the map $\Theta_r^*$. In particular, $P_0$ is the spinor space $P_0= \dneun$ and $P_1$ is the kernel of the Clifford multiplication $\Lambda^1 (\rneun) \otimes \dneun \to \dneun.$ The decomposition of the spaces $\Lambda^k (\rneun) \otimes \dneun$ is given by the formulas  (see \cite{10x})
$$ \Lambda^k (\rneun) \otimes \dneun = \sum\limits_{0 \le r \le \min (k,9-k)} \Theta_{k-1} \circ \ldots \circ \Theta_r (P_r) . $$

Moreover, $\Theta_{k-1} \circ \ldots \circ \Theta_{r} : P_r \to \Theta_{k-1} \circ \ldots \circ \Theta_r (P_r)$ is an isomorphism of $P_r$ onto the image. We apply this decomposition in the cases of  $k=1,2,3$:
\begin{enumerate}
\item $\Lambda^1 (\rneun) \otimes \dneun$ splits as a $\sneun$-representation into
$$ \Lambda^1 (\rneun) \otimes \dneun= \Theta_0 (P_0) \oplus P_1 = P_0 \oplus P_1 = \dneun \oplus P_1 ; $$
\item $\Lambda^2 (\rneun) \otimes \dneun$ splits as a $\sneun$-representation into
$$ \Lambda^2 (\rneun) \otimes \dneun = \Theta_1 \Theta_0 (P_0) \oplus \Theta_1 (P_1) \oplus P_2= \dneun \oplus P_1 \oplus P_2 ;  $$
\item $\Lambda^3 (\rneun) \otimes \dneun$ splits into
$$ \Lambda^3 (\rneun) \otimes \dneun = \dneun \oplus P_1 \oplus P_2 \oplus P_3 .$$
\end{enumerate}

The dimensions of the representations $P_r$ are given by:
$$ \dim P_0 = \dim \dneun = 16 \quad , \quad  \dim P_1 = 128 \, \, , \quad \dim P_2 = 432 \, \, , \quad \dim P_3 = 768.
$$

We decompose now the space $\Lambda^2 (\dneun)= \Lambda^2 ({\Bbb R}^{16})=$ {\mbox{\graf so}(16)} as well as  $\Lambda^3 (\Delta_9)$  into irreducible $\sneun$-components. It turns out that  only two components occur and these decompositions can be obtained in an elementary way.\\

{\bf Proposition 1:} {\it Under the action of the group $\sneun$ the spaces $\Lambda^p (\dneun)$ ($p=2,3$) decompose  into two irreducible components:}
$$ \Lambda^2 (\dneun) = \Lambda^2 ({\Bbb R}^9) \oplus \Lambda^3 (\rneun)  \quad , \quad \Lambda^3 (\Delta_9) = P_1 \oplus P_2 . $$

{\bf Proof:} We define an equivariant injection
$$ \Lambda^2 (\rneun) \longrightarrow \Lambda^2 (\dneun) $$

by the formula
$$ \sum\limits_{1 \le \alpha_1 < \alpha_2 \le 9} a_{\alpha_1 \alpha_2} \, \,  I_{\alpha_1} \wedge I_{\alpha_2} \longrightarrow \sum\limits_{1 \le \alpha_1 < \alpha_2 \le 9} a_{\alpha_1 \alpha_2} \, \, I_{\alpha_1}  I_{\alpha_2} . $$

Since $(I_{\alpha_1}  I_{\alpha_2})^* = I_{\alpha_2}^*  I_{\alpha_1}^* = I_{\alpha_2}  I_{\alpha_1} = - I_{\alpha_1}  I_{\alpha_2}$,  this map has values in the space $\Lambda^2 (\dneun)$ of all antisymmetric endomorphisms of $\dneun$. In a similar way we can define an injection $\Lambda^3 (\rneun) \to \Lambda^2 (\dneun)$ and, consequently, we have decomposed $\Lambda^2 (\dneun)$. $\Lambda^3 (\dneun)$ is the surjective image of the space
$$ \dneun \otimes \Lambda^2 (\dneun)= \dneun \otimes \Lambda^2 ({\Bbb R}^9) \oplus \dneun \otimes \Lambda^3 (\rneun) = 2 P_0 \oplus 2 P_1 \oplus 2 P_2 \oplus P_3 . $$

Therefore, $\Lambda^3 (\dneun)$ is a combination
$$ \Lambda^3 (\dneun)= A \cdot P_0 + B \cdot P_1 + C \cdot P_2 + D \cdot P_3 $$

where the integers $A,B,C \le 2$ are bounded by two. Inserting the dimensions  of the representations we immediately obtain that $(A, B, C, D)=(0,1,1,0)$ is the only possible solution.  \\
\mbox{}  \hfill \rule{3mm}{3mm}\\

The decomposition of the $\sneun$-representations $\Lambda^p (\dneun)$ for $p \ge 4$ is much more complicated and has been computed  by I. Agricola (see \cite{agricola}). We will use this result in the next section in an essential way.\\

{\bf Theorem 1 (see \cite{agricola}):} {\it The $\sneun$-representations $\Lambda^p (\dneun)$ are multiplicity-free. In particular, the representation $\Lambda^7 (\dneun)$ decomposes into
$$ \Lambda^7 (\dneun) = \dneun \oplus P_1 + P_2 + P_3 + \ldots  . $$}

{\bf Corollary 1:} {\it Let $L: \dneun \to \Lambda^7 (\dneun)$ be a linear $\sneun$-equivariant map. Then there exists a constant $C$ such that
$$ L(X) =C(X \, \mbox{}_{\_\!\!\_\!}{\scriptstyle{\mid}} \, \Omega_0^8)$$

holds for any $X \in \dneun$, where $\Omega_0^8 \in \Lambda^8 (\dneun)$ is the unique $\sneun$-invariant 8-form on $\dneun$.}\\

If the 16-dimensional manifold $\ms$ admits a $\sneun$-structure ${\cal R} \subset {\cal F} (\ms)$ with associated real vector bundle $\vn$, we will denote by ${\cal P}_r (\vn)$ the associated vector bundles
$$ {\cal P}_r (\vn)= {\cal R} \times_{\sneun} P_r . $$

The following bundle isomorphisms are consequences of the decompositions of the $\sneun$-representations.\\

{\bf Proposition 2:} {\it Let $\ms$ be a 16-dimensional Riemannian manifold with a fixed $\sneun$-structure. The following bundles are isomorphic:
\begin{enumerate}
\item $\Lambda^1 (\ms) \otimes \Lambda^1 (\vn) = \Lambda^1 (\ms) \oplus {\cal P}_1 (\vn)$; 
\item $\Lambda^1 (\ms) \otimes \Lambda^2 (\vn) = \Lambda^1 (\ms) \oplus {\cal P}_1 (\vn) \oplus {\cal P}_2 (\vn)$; 
\item $\Lambda^1 (\ms) \otimes \Lambda^3 (\vn) = \Lambda^1 (\ms) \oplus {\cal P}_1 (\vn) \oplus {\cal P}_2 (\vn) \oplus {\cal P}_3 (\vn)$ ; 
\item $\Lambda^2 (\ms) = \Lambda^2 (\vn) \oplus \Lambda^3 (\vn)$;
\item $\Lambda^3 (\ms) = {\cal P}_1 (\vn) \oplus {\cal P}_2 (\vn)$.
\end{enumerate}}

In a similar way we can study the 16-dimensional spin representations $\Delta_{16}^{\pm}$ as $\sneun$-representations. The element $e_1 \cdot \ldots e_{16}$ belongs to the subgroup $\widetilde{\sneun} \subset Spin(16)$ and acts on $\Delta_{16}^{\pm}$ by multiplication by $(\pm 1)$. Consequently, $\Delta_{16}^+$ is a $SO(9)$-representation, but $\Delta_{16}^-$ is a $\sneun$-representation. Both representations do not contain  non-trivial elements that are invariant under the $\sneun$-action. Indeed, a $\sneun$-invariant element in $\Delta_{16}^{\pm}$ would define a parallel spinor on the Cayley plane $F_4 / \sneun$. Since the Ricci tensor of this space is not zero, the spinor must vanish. Using these properties of the $\sneun$-representations $\Delta_{16}^{\pm}$ it is not hard to check that:
\begin{enumerate}
\item $\Delta_{16}^+$ is, as a $\sneun$-representation, isomorphic to $S_0^2 (\rneun) \oplus \Lambda^3 ({\Bbb R}^9)$, where $S_0^2$ denotes the space of traceless symmetric 2-tensors on $\rneun$.
\item $\Delta_{16}^-$ is the unique irreducible $\sneun$-representation of dimension 128. 
\end{enumerate}

This discussion yields the \\

{\bf Proposition 3:} {\it Let $\ms$ be a 16-dimensional Riemannian manifold with a fixed $\sneun$-structure. Then the spinor bundle $S^+ (\ms)$ of $\ms$ is isomorphic to
$$ S^+ (\ms) = S_0^2 (\vn) \oplus \Lambda^3 (\vn) . $$}

\section{The geometry of $\sneun$-structures}

In this section we introduce the 16 different types of weak geometric $\sneun$-structures. For this purpose, we will briefly recall how the different geometric classes of weak $G_2$-structures arise (weak geometric $U(n)$-structures can be treated in a completely similar way). Then we reformulate this scheme in a purely bundle theoretic way and use this approach to motivate the geometric types of weak $\sneun$-structures.\\

So, consider a 7-dimensional Riemannian manifold $(M^7 ,g)$ and a 3-form $\omega^3$ of general  type. This form defines a $G_2$-reduction ${\cal R}$ of its frame bundle ${\cal F} (M^7)$. Furthermore, its covariant derivative $\nabla \omega^3$ is a section in the bundle $T^* (M^7) \otimes \Lambda^3 (M^7)$ with special symmetry properties. Under the $G_2$-action, this space splits into 4 irreducible components, thus leading to 16 different geometric types of weak $G_2$-structures (see \cite{FG1}). Now we change the point of view. Denote by  $Z:T({\cal F}(M^7)) \to \mbox{{\graf so}(7)}$ the Levi-Civita connection of $M^7$. We decompose the Lie algebra into
$$ \mbox{{\graf so}(7)} = \mbox{\graf g}_2 \oplus \m , $$

when $\m \cong {\Bbb R}^7$ is the orthogonal complement of $\mbox{\graf g}_2$ inside $\mbox{{\graf so}(7)}$. The representation of $G_2$ on $\m$ is just the 7-dimensional standard representation of $G_2$. Given a $G_2$-structure defined by a subbundle ${\cal R} \subset {\cal F} (M^7)$,  we restrict $Z$ onto ${\cal R}$ and decompose it:
$$ Z_{|{\cal R}}= Z^* \oplus \Gamma . $$

Thus $Z^* :T({\cal R}) \to \mbox{\graf g}_2$ is a connection in the $G_2$-principal fibre bundle ${\cal R}$ and $\Gamma$ is a 1-form on $M^7$ with values in the associated vector bundle
$$ {\cal R} \times_{G_2} \m = {\cal R} \times_{G_2} {\Bbb R}^7 = T(M^7) . $$

Denote by $\rho_3$ the representation of $SO(7)$ on $\Lambda^3 ({\Bbb R}^7)$. Then the covariant derivative $\nabla \omega^3$ is given by
$$ \nabla \omega^3 = \rho_3 (\Gamma)(\omega^3) . $$

Consider now the $G_2$-equivariant map
$$ {\Bbb R}^7 \otimes {\Bbb R}^7 \to {\Bbb R}^7 \otimes \Lambda^3 ({\Bbb R}^3) \, \, , \quad X \otimes Y \mapsto X \otimes (\rho_3 (Y) \omega^3) . $$

The crucial point is that this map is injective. Consequently, the $G_2$-type of $\nabla \omega^3$ is uniquely determined by the $G_2$-type of $\Gamma \in \Lambda^1 (M^7) \otimes T(M^7)$. Since ${\Bbb R}^7 \otimes {\Bbb R}^7$ splits again into 4 summands, we reobtain the previous 16 classes of weak geometric $G_2$-structures, and may view this as an alternative, but completely equivalent definition of these structures. The advantage of this approach is that the form $\omega^3$ does not appear any more in the definition and can thus be used, in our situation, for defining weak $\sneun$-structures in 16 dimensions.\\ 

Let $(\ms, g)$ be an oriented, 16-dimensional Riemannian manifold and fix some $\sneun$-structure ${\cal R} \subset {\cal F} (\ms)$. We decompose the Lie algebra $\so$ into 
$$ \so = {\mbox{\graf spin}(9)} \oplus \m . $$

Using the symmetric operators $I_{\alpha} : {\Bbb R}^{16} \to {\Bbb R}^{16}$ we know that $I_{\alpha}  I_{\beta}$ and $I_{\alpha}  I_{\beta}  I_{\gamma}$ \, $(\alpha < \beta < \gamma)$ are antisymmetric. The Lie algebra $\mbox{\graf spin}(9)$ is spanned by the elements $I_{\alpha}  I_{\beta} \in \Lambda^2 ({\Bbb R}^{16})$ and the operators $I_{\alpha}  I_{\beta}  I_{\gamma}$ form a basis of the orthogonal complement $\m$:
\begin{eqnarray*}
 {\mbox{\graf spin}} (9) &=& {\mathrm{Lin}} (I_{\alpha}  I_{\beta} :\quad \alpha < \beta)\\
 \m &=& \mbox{$\mathrm{Lin}$} (I_{\alpha}  I_{\beta}  I_{\gamma} : \quad \alpha < \beta < \gamma) . 
\end{eqnarray*}

The $\sneun$-representation $\mbox{Ad}$ in the complement $\m$ is equivalent to the representation $\Lambda^2 ({\Bbb R}^9)$.  The Levi-Civita connection of the Riemannian manifold is a 1-form on ${\cal F} (\ms)$
$$ Z:T({\cal F} (\ms)) \to \so $$

with values in $\so$. We restrict the connection form to the $\sneun$-structure and decompose it with respect to the decomposition of the Lie algebra $\so$ into
$$ Z_{|{\cal R}}= Z^* \oplus \Gamma  . $$

Then, $Z^*$ is a connection in the principal $\sneun$-fibre bundle ${\cal R}$,  and $\Gamma$ is a tensorial 1-form of type Ad. Therefore, $\Gamma$ is a 1-form defined on $\ms$ with values in the associated bundle 
$$ {\cal R} \times_{\sneun} \m ={\cal R} \times_{\sneun} \Lambda^3 ({\Bbb R}^9)= \Lambda^3 (\vn) . $$

In case $\Gamma \equiv 0$,  the Riemannian manifold has a holonomy group contained in $\sneun$ and a classical result of Alekseevski/Brown/Gray (see \cite{11x}, \cite{x}) states that $\ms$ is either flat or isometric to one of the symmetric spaces $F_4 / \sneun$ or $F_4^* / \sneun$. On the other hand, $\Lambda^1 (\ms) \otimes \Lambda^3 (\vn)$ splits into 4 subbundles (see Section 7). From this point of view we obtain 16 classes of ''weak $\sneun$-structures'' depending on the algebraic type of $\Gamma$. In this paper we will only study one class of these $\sneun$-structures. The splitting 
$$ \Lambda^1 (\ms) \otimes \Lambda^3 (\vn) = \Lambda^1 (\ms) \oplus {\cal P}_1 (\vn) \oplus {\cal P}_2 (\vn) \oplus {\cal P}_3 (\vn) $$

is the background for the following\\

{\bf Definition:} A $\sneun$-structure on a 16-dimensional Riemannian manifold is called {\it nearly parallel} if $\Gamma$ is a vector field.\\

The inclusion $\Lambda^1 (\ms) \subset \Lambda^1 (\ms) \otimes \Lambda^3 (\vn) \subset \Lambda^1 (\ms) \otimes \Lambda^2 (\ms)$ is given by the formula
$$ \Gamma \longrightarrow 6 \sum\limits_{1 \le \alpha < \beta < \gamma \le 9} I_{\gamma} I_{\beta} I_{\alpha} (\Gamma) \otimes (I_{\alpha}  I_{\beta}  I_{\gamma})$$

and, therefore, in case of a nearly parallel $\sneun$-structure this sum coincides with the difference $Z - Z^*$ of the two connections. Since the Levi-Civita connection is a torsion free connection, we obtain the following formula for the torsion tensor $T^{\star}$ of the connection $Z^{\star}$
$$ T^{\star} (X,Y) = 6 \sum\limits_{1 \le \alpha < \beta < \gamma \le 9} \Big\{ g(\Gamma, I_{\alpha} I_{\beta} I_{\gamma} (X)) I_{\alpha} I_{\beta} I_{\gamma} (Y) - g(\Gamma, I_{\alpha} I_{\beta} I_{\gamma} (Y)) I_{\alpha} I_{\beta} I_{\gamma} (X) \Big\} . $$

$X, Y \in T(\ms)$ are vectors tangent to $\ms$. In particular, we have
$$ g(T^* (X,Y), \Gamma) \equiv 0 . $$

Fix a local section $(e_1, \ldots , e_{16})$ in the reduction ${\cal R} \subset {\cal F} (\ms)$ and denote by $\sigma^1 , \ldots , \sigma^{16}$ the dual frame. Then the 1-form 
$$ \sum\limits_{1 \le \alpha < \beta < \gamma \le 9} I_{\gamma} I_{\beta} I_{\alpha} (\Gamma) \otimes (I_{\alpha}  I_{\beta}  I_{\gamma})$$

with values in $\Lambda^2 (\ms)$ is given by the $(16\times 16)$-matrix $(\Gamma= e_{16})$:\\

\setlength{\tabcolsep}{0.9mm}
{\tiny
\begin{tabular}{|c|c|c|c|c|c|c|c|c|c|c|c|c|c|c|c|}
\hline &&&&&&&&&&&&&&&\\ 
0 & $2 \sigma^{15}$ &  $-2 \sigma^{14}$ & $-2 \sigma^{13}$ & $2 \sigma^{12}$ & $2 \sigma^{11}$ & $- 2 \sigma^{10}$ & $ -2 \sigma^9$ & $- \sigma^8$ & $- \sigma^7$ & $ \sigma^6$ & $\sigma^5$ & $- \sigma^4$ & $ - \sigma^3$ & $\sigma^2$ & $7 \sigma^1$\\ &&&&&&&&&&&&&&&\\
\hline &&&&&&&&&&&&&&&\\
$-2 \sigma^{15}$ & 0 & $2 \sigma^{13}$ & $ -2 \sigma^{14}$ & $-2 \sigma^{11}$ & $2 \sigma^{12}$ & $2 \sigma^9$ & $-2 \sigma^{10}$ & $\sigma^7$ & $- \sigma^8$ & $- \sigma^5$ & $\sigma^6$ & $\sigma^3$ & $- \sigma^4$ & $- \sigma^1$ & $7 \sigma^2$\\ &&&&&&&&&&&&&&&\\ \hline
&&&&&&&&&&&&&&&\\
$2 \sigma^{14}$ & $- \sigma^{13}$ & $0$ & $ -2 \sigma^{15}$ & $2 \sigma^{10}$ & $-2 \sigma^{9}$ & $2 \sigma^{12}$ & $-2 \sigma^{11}$ & $-\sigma^6$ & $ \sigma^5$ & $- \sigma^8$ & $\sigma^7$ & $-\sigma^2$ & $\sigma^1$ & $- \sigma^4$ & $7 \sigma^3$\\ &&&&&&&&&&&&&&&\\ \hline
&&&&&&&&&&&&&&&\\
$2 \sigma^{13}$ & $2 \sigma^{14}$ & $2 \sigma^{15}$ & $0$ & $-2 \sigma^{9}$ & $-2 \sigma^{10}$ & $-2 \sigma^{11}$ & $-2 \sigma^{12}$ & $-\sigma^5$ & $- \sigma^6$ & $- \sigma^7$ & $-\sigma^8$ & $\sigma^1$ & $\sigma^2$ & $\sigma^3$ & $7 \sigma^4$\\ &&&&&&&&&&&&&&&\\ \hline
&&&&&&&&&&&&&&&\\
$-2 \sigma^{12}$ & $2 \sigma^{11}$ & $-2 \sigma^{10}$ & $2 \sigma^{9}$ & $0$ & $-2 \sigma^{15}$ & $2 \sigma^{14}$ & $-2 \sigma^{13}$ & $\sigma^4$ & $- \sigma^3$ & $\sigma^2$ & $-\sigma^1$ & $-\sigma^8$ & $\sigma^7$ & $- \sigma^6$ & $7 \sigma^5$\\ &&&&&&&&&&&&&&&\\ \hline
&&&&&&&&&&&&&&&\\
$-2 \sigma^{11}$ & $-2 \sigma^{12}$ & $2 \sigma^{9}$ & $2 \sigma^{10}$ & $2 \sigma^{15}$ & 0 & $-2 \sigma^{13}$ & $-2 \sigma^{14}$ & $\sigma^3$ & $\sigma^4$ & $- \sigma^1$ & $-\sigma^2$ & $-\sigma^7$ & $- \sigma^8$ & $\sigma^5$ & $7 \sigma^6$\\ &&&&&&&&&&&&&&&\\ \hline
&&&&&&&&&&&&&&&\\
$2 \sigma^{10}$ & $-2 \sigma^9$ & $-2 \sigma^{12}$ & $2 \sigma^{11}$ & $-2\sigma^{14}$ & $\sigma^{13}$ & $0$ & $-2 \sigma^{15}$ & $-\sigma^2$ & $\sigma^1$ & $\sigma^4$ & $-\sigma^3$ & $\sigma^6$ & $- \sigma^5$ & $- \sigma^8$ & $7 \sigma^7$\\ &&&&&&&&&&&&&&&\\ \hline
&&&&&&&&&&&&&&&\\
$2 \sigma^{9}$ & $2 \sigma^{10}$ & $2 \sigma^{11}$ & $ 2 \sigma^{12}$ & $2 \sigma^{13}$ & $2 \sigma^{14}$ & $2 \sigma^{15}$ & $0$ & $\sigma^1$ & $\sigma^2$ & $\sigma^3$ & $\sigma^4$ & $\sigma^5$ & $\sigma^6$ & $\sigma^7$ & $7 \sigma^8$\\ &&&&&&&&&&&&&&&\\ \hline
&&&&&&&&&&&&&&&\\
$\sigma^{8}$ & $- \sigma^7$ & $\sigma^{6}$ & $\sigma^{5}$ & $- \sigma^{4}$ & $-\sigma^{3}$ & $\sigma^2$ & $- \sigma^{1}$ & $0$ & $0$ & $0$ & $0$ & $0$ & $0$ & $0$ & $4 \sigma^9$\\ &&&&&&&&&&&&&&&\\ \hline &&&&&&&&&&&&&&&\\
$\sigma^{7}$ & $\sigma^8$ & $-\sigma^{5}$ & $\sigma^{6}$ & $\sigma^{3}$ & $-\sigma^{4}$ & $-\sigma^1$ & $- \sigma^{2}$ & $0$ & $0$ & $0$ & $0$ & $0$ & $0$ & $0$ & $4 \sigma^{10}$\\ &&&&&&&&&&&&&&&\\ \hline &&&&&&&&&&&&&&&\\
$-\sigma^{6}$ & $\sigma^5$ & $\sigma^{8}$ & $\sigma^{7}$ & $- \sigma^{2}$ & $\sigma^{1}$ & $-\sigma^4$ & $- \sigma^{3}$ & $0$ & $0$ & $0$ & $0$ & $0$ & $0$ & $0$ & $4 \sigma^{11}$\\ &&&&&&&&&&&&&&&\\ \hline &&&&&&&&&&&&&&&\\
$-\sigma^{5}$ & $- \sigma^6$ & $-\sigma^{7}$ & $\sigma^{8}$ & $\sigma^{1}$ & $\sigma^{2}$ & $\sigma^3$ & $- \sigma^{4}$ & $0$ & $0$ & $0$ & $0$ & $0$ & $0$ & $0$ & $4 \sigma^{12}$\\ &&&&&&&&&&&&&&&\\ \hline &&&&&&&&&&&&&&&\\
$\sigma^{4}$ & $- \sigma^3$ & $\sigma^{2}$ & $-\sigma^{1}$ & $\sigma^{8}$ & $\sigma^{7}$ & $-\sigma^6$ & $- \sigma^{5}$ & $0$ & $0$ & $0$ & $0$ & $0$ & $0$ & $0$ & $4 \sigma^{13}$\\ &&&&&&&&&&&&&&&\\ \hline &&&&&&&&&&&&&&&\\
$\sigma^{3}$ & $\sigma^4$ & $-\sigma^{1}$ & $-\sigma^{2}$ & $- \sigma^{7}$ & $\sigma^{8}$ & $\sigma^5$ & $- \sigma^{6}$ & $0$ & $0$ & $0$ & $0$ & $0$ & $0$ & $0$ & $4 \sigma^{14}$\\ &&&&&&&&&&&&&&&\\ \hline &&&&&&&&&&&&&&&\\
$-\sigma^{2}$ & $\sigma^1$ & $\sigma^{4}$ & $-\sigma^{3}$ & $\sigma^{6}$ & $-\sigma^{5}$ & $\sigma^8$ & $- \sigma^{7}$ & $0$ & $0$ & $0$ & $0$ & $0$ & $0$ & $0$ & $4 \sigma^{15}$\\ &&&&&&&&&&&&&&&\\ \hline &&&&&&&&&&&&&&&\\
$-7 \sigma^{1}$ & $-7 \sigma^{2}$ & $-7 \sigma^{3}$ & $ -7 \sigma^{4}$ & $-7 \sigma^{5}$ & $-7 \sigma^{6}$ & $-7 \sigma^{7}$ & $-7\sigma^{8}$ & $-4\sigma^9$ & $-4\sigma^{10}$ & $-4\sigma^{11}$ & $-4\sigma^{12}$ & $-4\sigma^{13}$ & $-4\sigma^{14}$ & $-4\sigma^{15}$ & $0$\\ &&&&&&&&&&&&&&&\\ \hline
\end{tabular}
}

\vspace{0.5cm}

Let us explain the representation of the endomorphisms $I_{\alpha}$ we used here. Denote by $E_{ij} \in \, \mbox{\graf so(8)}$ the standard basis of the Lie algebra $\mbox{\graf so(8)}$ and consider the real representation of the 7-dimensional Clifford algebra (see \cite{BFGK}):
\begin{eqnarray*}
e_1 &=& \hspace{0.3cm} E_{18} + E_{27} - E_{36} - E_{45}\\
e_2 &=& -E_{17} + E_{28} + E_{35} - E_{46}\\
e_3 &=& -E_{16} + E_{25} - E_{38} + E_{47}\\
e_4 &=& -E_{15} - E_{26} - E_{37} - E_{48}\\
e_5 &=& - E_{13} - E_{24} + E_{57} + E_{68}\\
e_6 &=& \hspace{0.3cm} E_{14}  - E_{23} - E_{58} + E_{67}\\
e_7 &=& \hspace{0.3cm} E_{12} - E_{34} - E_{56} + E_{78} \quad . \\
\end{eqnarray*}

A representation of the matrices $I_1 , \ldots , I_9$ in ${\Bbb R}^{16}$ is, for example, given by
$$ I_{\alpha} = \left( \begin{array}{cc} 0&-e_{\alpha}\\ e_{\alpha} &0 \end{array} \right) \quad 1 \le \alpha \le 7 \quad , \quad  I_{8} = \left( \begin{array}{cc} 0&E\\ E &0 \end{array} \right) \quad  , \quad  I_{9} = \left( \begin{array}{cc} E&0\\ 0 &-E \end{array} \right) , $$

where $E$ is the identity on ${\Bbb R}^8$. A computer calculation immediately yields the formula for the $(16 \times 16)$-matrix.\\

{\bf Example:} The manifold $\ms = S^{15} \times S^{1}$ admits a homogeneous, nearly parallel $\sneun$-structure such that $0 \not= \Gamma \in \Lambda^1 (\ms)$. Indeed, $\ms$ is a homogeneous space
$$ \ms = (\sneun / Spin_{\Delta} (7)) \times S^{1} =   (SO(16)/SO(15)) \times S^1. $$

The frame bundle ${\cal F} (\ms)$ admits the $SO(15)$-reduction ${\cal F}_0 (\ms) =   SO(16) \times S^1 \subset {\cal F} (\ms)$. Denote by $\Theta$ the Maurer-Cartan form of the Lie group $SO(16)$. The Levi-Civita connection of $\ms$ is given by the 1-form
$$ Z= pr_{\mbox{\graf so}(15)} (\Theta) : T( SO(16) \times S^1) \to {\mbox{\graf so}(15)} , $$

where we decompose ${\mbox{\graf so}(16)}$ into ${\mbox{\graf so}(15)}$ and its orthogonal complement
$$ {\mbox{\graf so}(16)} = {\mbox{\graf so}(15)} \oplus \n  $$

and project onto ${\mbox{\graf so}(15)}$. On the other hand, we have the commutative diagram
$$
\begin{diagram}
\node{\sneun} \arrow{e} \node{SO(16)} \\
\node{Spin_{\Delta}(7)} \arrow{n} \arrow{e} \node{SO(15)  . } \arrow{n}
\end{diagram} 
$$

${\cal F} (\ms)$ admits the $Spin_{\Delta} (7)$-reduction ${\cal R} =  \sneun \times S^1$. We restrict the Levi-Civita connection to ${\cal R}$ and decompose it with respect to the orthogonal decomposition
$$ {\mbox{\graf so}(16)} = \mbox{\graf spin}(9) \oplus \m . $$

Finally, let us decompose the Lie algebra
$$ {\mbox{\graf spin}(9)} = \mbox{\graf spin}_{\Delta} (7) \oplus {\mbox{\graf k}} . $$

Using the $Spin_{\Delta} (7)$-reduction ${\cal R}$ the 1-form $\Gamma$ is given by the formula
$$\Gamma :T({\cal R}) \to \m \quad , \quad \Gamma = pr_{\m} \circ pr_{\mbox{\graf so}(15)} (\Theta) . $$

The map $pr_{\m} \circ pr_{\mbox{\graf so}(15)} : {\mbox{\graf spin}(9)} \to \m $ vanishes on the subspace $\mbox{\graf spin}_{\Delta}(7)$ and, consequently, it is a $Spin_{\Delta}(7)$-equivariant map $\mbox{\graf k} \to \m$. We compute the formula for $\Gamma$. Consider an element $W$ in the Lie algebra $\mbox{\graf spin}(9)$,
$$ W = \sum\limits_{1 \le \alpha < \beta \le 9} \mbox{x}_{\alpha \beta} I_{\alpha}  I_{\beta} $$

and let us introduce $\mu^1 , \ldots , \mu^{15}$ given by\\

\hspace{2cm}$ \mu^1 = 2 \, \mbox{x}_{19} \hspace{0.7cm}  , \hspace{0.25cm} \quad \mu^2 = 2 \, \mbox{x}_{29}  \quad \hspace{0.55cm}  , \quad  \hspace{0.5cm} \mu^3 = - 2 \, \mbox{x}_{39} \hspace{0.3cm}  , \hspace{0.5cm} \quad \mu^4 = - 2 \, \mbox{x}_{49}$\\

\hspace{2cm}$ \mu^5 = - 2 \, \mbox{x}_{69} \quad , \quad \hspace{0.3cm} \mu^6 = 2 \, \mbox{x}_{59} \quad \hspace{0.5cm} , \hspace{0.9cm} \mu^7 = 2 \, \mbox{x}_{79} \quad \hspace{0.2cm} , \quad \hspace{0.5cm} \mu^8 = - 2 \, \mbox{x}_{89} , $\\

\hspace{2cm}$ \mu^9 = 2 \, \mbox{x}_{18} + 2 \, \mbox{x}_{27} + 2 \, \mbox{x}_{35} - 2 \, \mbox{x}_{46} \quad \hspace{0.2cm} , \quad \hspace{0.5cm}    \mu^{10} = - 2 \, \mbox{x}_{17} + 2 \, \mbox{x}_{28} + 2 \, \mbox{x}_{36} + 2 \, \mbox{x}_{45} $\\

\hspace{2cm}$ \mu^{11} = 2 \, \mbox{x}_{15} + 2 \, \mbox{x}_{26} - 2 \, \mbox{x}_{38} + 2 \, \mbox{x}_{47} \quad  , \quad \mu^{12} = - 2 \, \mbox{x}_{16} + 2 \, \mbox{x}_{25} - 2 \, \mbox{x}_{37} - 2 \, \mbox{x}_{48} $\\

\hspace{2cm}$ \mu^{13} = 2 \, \mbox{x}_{14} - 2 \, \mbox{x}_{23} - 2 \, \mbox{x}_{57} - 2 \, \mbox{x}_{68} \quad , \quad \mu^{14} = 2 \, \mbox{x}_{13} + 2 \, \mbox{x}_{24} + 2 \, \mbox{x}_{58} - 2 \, \mbox{x}_{67} $\\

\hspace{2cm}$ \mu^{15} = 2 \, \mbox{x}_{12} - 2 \, \mbox{x}_{34} - 2 \, \mbox{x}_{56} + 2 \, \mbox{x}_{78}.$\\

Then a direct calculation yields that $\Gamma (W)$ is given by the same $(16 \times 16)$-matrix as the endomorphism $\sum\limits_{1 \le \alpha < \beta <€\gamma \le 9} I_{\gamma} I_{\beta} I_{\alpha} (e_{16}) \otimes I_{\alpha}  I_{\beta}  I_{\gamma}$ if we replace the form $\sigma^i$ by $\mu^i$. \\

Consequently, the $\sneun$-structure on $S^1 \times S^{15}$ is nearly parallel and the vector field $\Gamma = \frac{\partial}{\partial \varphi}$ is the unit vector field tangent to $S^1$. \hfill \rule{3mm}{3mm}\\

\newcommand{\k}{\mbox{\graf k}}

\vspace{0.5cm}

We derive now a differential equation for the canonical 8-form $\Omega^8$ of a nearly parallel $\sneun$-structure. This equation - in contrast to the case of weak $G_2$-structures (see \cite{FG1}) or weak $Spin(7)$-structures (see \cite{F2}) - does not completely characterize the nearly parallel $\sneun$-structures. We will contract the covariant derivative $\nabla \Omega^8 \in \Gamma (\Lambda^1 (\ms) \otimes \Lambda^8 (\ms))$ once in order to obtain an equation for $\delta \Omega^8$.\\

{\bf Theorem 1:} {\it Let ${\cal R} \subset {\cal F} (\ms)$ be a nearly parallel $\sneun$-structure on $\ms$ with vector field $\Gamma$. Then the canonical form $\Omega^8$ satisfies the equations
$$ \delta \Omega^8 = - 504 (\Gamma \, \mbox{}_{\_\!\!\_\!}{\scriptstyle{\mid}} \, \Omega^8) \quad , \quad d \Omega^8 = - 504 \star (\Gamma \, \mbox{}_{\_\!\!\_\!}{\scriptstyle{\mid}} \, \Omega^8) . $$}

{\bf Proof:} The 8-form $\Omega^8$ is defined by a $\sneun$-invariant form in $\Lambda^8 (\dneun)= \Lambda^8 ({\Bbb R}^{16})$. Therefore, $\Omega^8$ is parallel with respect to the connection $Z^*$, 
$$ \nabla^* \Omega^8 = 0 . $$

The covariant derivative with respect to the Levi-Civita connection is now given by the formula
$$ \nabla \Omega^8 = 6 \sum\limits_{1 \le \alpha < \beta < \gamma \le 9} I_{\gamma} I_{\beta} I_{\alpha} (\Gamma) \otimes \rho_8 (I_{\alpha}  I_{\beta}  I_{\gamma} ) \Omega^8 ,  $$

where $\rho_8 : \mbox{\graf so}(16) \to \mbox{\graf so} (\Lambda^8 ({\Bbb R}^{16}))$ is the differential of the representation of the group $SO(16)$ in $\Lambda^8 ({\Bbb R}^{16})$. Contracting this equation we obtain that
$$ \delta \Omega^8 = - 6 \sum\limits_{1 \le \alpha < \beta < \gamma \le 9} I_{\gamma} I_{\beta} I_{\alpha} (\Gamma) \, \mbox{}_{\_\!\!\_\!}{\scriptstyle{\mid}} \, ( \rho_8 (I_{\alpha}  I_{\beta}  I_{\gamma} ) \Omega^8 ) $$

is defined by a linear, $\sneun$-equivariant map of ${\Bbb R}^{16} = \dneun$ into $\Lambda^7 (\dneun) = \Lambda^7 ({\Bbb R}^{16})$. Using Theorem 1 of Section 7 we conclude that there exists a constant $C$ such that 
$$ \delta \Omega^8 = C ( \Gamma \, \mbox{}_{\_\!\!\_\!}{\scriptstyle{\mid}} \, \Omega^8) $$

holds. Taking into account the explicit formula for $\rho_8$ we obtain $C= - 504$. \hfill \rule{3mm}{3mm}\\

{\bf Remark:} A second natural class of $\sneun$-structures occurs if $\Gamma$ belongs to the subbundle ${\cal P}_1 (V^9) \oplus {\cal P}_2 (V^9) \subset \Lambda^1 (M^{16}) \otimes \Lambda^3 (V^9)$. According to Proposition 2 of Section 7 $\Gamma$ is a pair of 3-forms $\Gamma = \Gamma_1 + \Gamma_2 \in {\cal P}_1 (V^9) \oplus {\cal P}_2 (V^9) = \Lambda^3 (M^{16})$. A $\sneun$-structure of this type is called a {\it nearly 3-parallel $\sneun$-structure}. Let us derive the corresponding differential equation for $d \Omega^8$. Remark that there are two non-trivial $\sneun$-equivariant linear maps of $\Lambda^3 (\Delta_9)$ into $\Lambda^9 (\Delta_9)$:
\begin{eqnarray*}
\Psi_1 (\Gamma) &=& \sum\limits^{16}_{i=1} (e_i \, \mbox{}_{\_\!\!\_\!}{\scriptstyle{\mid}} \,  \Gamma) \wedge (e_i \, \mbox{}_{\_\!\!\_\!}{\scriptstyle{\mid}} \,  \Omega^8)\\
\Psi_2 (\Gamma) &=& \star \Big( \sum\limits^{16}_{i,j=1} (e_i \, \mbox{}_{\_\!\!\_\!}{\scriptstyle{\mid}} \,  e_j \, \mbox{}_{\_\!\!\_\!}{\scriptstyle{\mid}} \,  \Gamma) \wedge (e_i \, \mbox{}_{\_\!\!\_\!}{\scriptstyle{\mid}} \,  e_j \, \mbox{}_{\_\!\!\_\!}{\scriptstyle{\mid}} \,  \Omega^8) \Big) . 
\end{eqnarray*}

Here $\star$ denotes the Hodge operator acting on forms in sixteen variables. Since $\Lambda^9 (\Delta_9) \simeq \Lambda^7 (\Delta_9)$ is a multiplicity-free $\sneun$-representation and the representation $\Lambda^3 (\Delta_9) = {\cal P}_1 ({\Bbb R}^9) \oplus {\cal P}_2 ({\Bbb R}^9)$ splits into two irreducible components, any $\sneun$-equivariant map $\Lambda^3 (\Delta_9) \to \Lambda^9 (\Delta_9)$ is a linear combination of $\Psi_1$ and $\Psi_2$. The parameters depend on the fixed isomorphism ${\cal P}_1 ({\Bbb R}^9) \oplus {\cal P}_2 ({\Bbb R}^9) = \Lambda^3 (\Delta_9)$ and they can be normalized to one. Therefore, we obtain two differential equations depending on the algebraic type of the 3-form $\Gamma$:
\begin{eqnarray*}
d \Omega^8 &=& \sum\limits^{16}_{i=1} (e_i \, \mbox{}_{\_\!\!\_\!}{\scriptstyle{\mid}} \,  \Gamma) \wedge (e_i \, \mbox{}_{\_\!\!\_\!}{\scriptstyle{\mid}} \,  \Omega^8) \\
d \Omega^8 &=& \star \Big( \sum\limits^{16}_{i,j=1} (e_i \, \mbox{}_{\_\!\!\_\!}{\scriptstyle{\mid}} \,  e_j \, \mbox{}_{\_\!\!\_\!}{\scriptstyle{\mid}} \,  \Gamma) \wedge (e_i \, \mbox{}_{\_\!\!\_\!}{\scriptstyle{\mid}} \,  e_j \, \mbox{}_{\_\!\!\_\!}{\scriptstyle{\mid}} \,  \Omega^8) \Big) . 
\end{eqnarray*}  

\section{The twistor space of a $\sneun$-structure}

Let ${\cal C}$ be the space of all complex structures ${\cal J}$ compatible with the metric and acting on ${\Bbb R}^{16} = \dneun$:
$$ {\cal J}^2 = - \mbox{Id} \, \quad , \quad {\cal J} {\cal J}^* = \mbox{Id} . $$

We consider the intersection of ${\cal C}$ with the space $\Lambda^2 ({\Bbb R}^9)= {\mbox{\graf spin}(9)}$:
$$ {\cal T}_1 = {\cal C} \cap \Lambda^2 ({\Bbb R}^9) = \Big\{ {\cal J} = \sum\limits_{1 \le \alpha < \beta \le 9} \mbox{x}_{\alpha \beta} I_{\alpha} I_{\beta} : \, \, {\cal J}^2 = - \mbox{Id} \Big\} . $$

{\bf Proposition 1:} {\it The group $\sneun$ acts transitively on ${\cal T}_1$. The isotropy group of the operator $I_1 I_2 \in {\cal T}_1$ is isomorphic to}
$$ \sneun \cap U(8) = {Spin} (2) \times_{{\Bbb Z}_2} {Spin} (7) . $$

{\bf Proof:} The normal form of a 2-form is
$$ {\cal J} = a I_1 I_2 + b I_3 I_4 + c I_5 I_6 + d I_7 I_8 . $$

The condition ${\cal J}^2= - \mbox{Id}$ yields the equations
$$ (a+b+c-d)^2 =1 \quad \quad \quad (a+b+c+d)^2 =1$$
$$(a+b-c+d)^2 =1 \quad \quad \quad (a+b-c-d)^2 =1 $$
$$ (a-b+c+d)^2 =1 \quad \quad \quad (a-b+c-d)^2 =1 $$
$$(-a+b+c+d)^2 =1 \quad \quad (-a+b+c-d)^2 =1 . $$

The solutions of these equations are the 4-tuples $(\pm 1,0,0,0), (0, \pm 1,0,0), (0,0, \pm 1,0)$ and $(0,0,0,\pm 1)$. Therefore, $\sneun$ acts transitively on ${\cal T}_1$. \mbox{} \hfill \rule{3mm}{3mm}\\

{\bf Corollary:} {\it The space ${\cal T}_1$ is a complex manifold isomorphic to the quadric $Q$ in ${\Bbb P}^8 ({\Bbb C})$ defined by the equation $z_0^2 + \ldots + z_8^2=0$:}
$$ {\cal T}_1 = \sneun / Spin(2) \times_{{\Bbb Z}_2} Spin(7) = Q . $$

{\bf Definition:} {\it Let $(M^{16},g)$ be an oriented, 16-dimensional Riemannian manifold with a fixed $\sneun$-structure ${\cal R}  \subset {\cal F} (M^{16})$. We define the twistor space ${\cal T}_1 (M^{16})$ as the associated bundle}
$$ {\cal T}_1 (M^{16}) = R \times_{Spin(9)} {\cal T}_1 . $$

The general twistor construction (see \cite{BBO}, \cite{OBR}) yields a natural almost complex structure ${\cal J}$ on the twistor space ${\cal T}_1 (M^{16})$. ${\cal J}$ depends on the restriction $Z^*$ of the Levi-Civita connection $Z$. There are two types of integrability conditions for ${\cal J}$. The first type is an algebraic condition for the torsion tensor $T^*$ of the connection $Z^*$. The identity \\

$(\star)$ \hfill $T^* ({\cal J} X , {\cal J}Y) - {\cal J} T^* ({\cal J} X,Y) - {\cal J}T^* (X, {\cal J} Y) - T^*(X,Y) =0 $ \hfill \mbox{}\\

should be satisfied for any operator ${\cal J} \in {\cal C} \cap \Lambda^2 (V^9)$ and any pair of vectors $X,Y$.\\

{\bf Proposition 2:} {\it Let ${\cal R} \subset {\cal F}(M^{16})$ be a nearly parallel $\sneun$-structure on a 16-dimensional Riemannian manifold. Then the torsion tensor $T^*$ satisfies the integrability condition $(\star)$.}\\

{\bf Proof:} A direct calculation using the explicit formula of the torsion tensor $T^*$ in case $\Gamma$ is a vector yields the result. \mbox{} \hfill \rule{3mm}{3mm}\\

The second integrability condition is an identity for the curvature $R^*$ of the connection $Z^*$:
$$ (\star \star) \quad \quad [ R^* ({\cal J} X, {\cal J}Y),  {\cal J}] - {\cal J} [R^*({\cal J}X,Y), {\cal J}]  - {\cal J} [R^* (X, {\cal J}X) , {\cal J}] - [R^* (X,Y) , {\cal J}]=0 $$

for any pair of vectors $X,Y$ and for any operator ${\cal J} \in {\cal C} \cap \Lambda^2 (V^9)$. The relation $Z=Z^* + \Gamma$ yields the following formula for the curvature forms $\Omega^{Z^*}$ and $\Omega^Z$ of the connections $Z^*$ and $Z$:
$$ \Omega^{Z^*} = \Omega^Z - D^Z (\Gamma) + \frac{1}{2} [\Gamma , \Gamma] . $$

Here $D^Z (\Gamma)=d \Gamma + [Z , \Gamma]$ denotes the derivative of the 1-form $\Gamma$ with respect to the Levi-Civita connection.\\

{\bf Proposition 3:} {\it Let ${\cal R} \subset {\cal F} (M^{16})$ be a nearly parallel $\sneun$-structure on a 16-dimensional Riemannian manifold. Then the 2-form $[\Gamma , \Gamma]$ satisfies the condition $(\star \star)$ for any ${\cal J} \in {\cal C} \cap \Lambda^2 (V^9)$. }\\

{\bf Proof:} A direct calculation yields the result. \hfill \rule{3mm}{3mm}\\

Let us introduce a 2-form $W$ with values in the bundle of endomorphisms of $T(M^{16})$
$$ W= \Omega^Z - D^Z (\Gamma) . $$

Then we obtain the following \\

{\bf Theorem 1:}€\, {\it The twistor space ${\cal T}_1 (M^{16})$ of a 16-dimensional Riemannian manifold with nearly parallel $\sneun$-structure is a complex manifold if and only if the tensor field $W= \Omega^Z - D^Z (\Gamma)$ satisfies the equation 
$$ [W ({\cal J}X,{\cal J}Y) , {\cal J}] - {\cal J}[W({\cal J}X,Y) , {\cal J}] -  {\cal J} [W(X, {\cal J}Y) , {\cal J}] - [W(X,Y) , {\cal J}]=0 $$

for any point ${\cal J} \in {\cal T}_1 (M^{16})$ of the twistor space.}\\

We formulated the integrability condition for the almost complex structure of ${\cal T}_1 (M^{16})$ viewing the tensor  $W= \Omega^Z - D^Z (\Gamma)$ as a $(2,0)$-tensor with values in the bundle of antisymmetric endomorphisms of $T(M^{16})$. Remark that $W$ does not satisfy the first and second Bianchi identity in general. We can understand $W$ as an endomorphism acting on the bundle of 2-forms,
$$ W : \Lambda^2 (M^{16}) \to \Lambda^2 (M^{16}) . $$

In this case the integrability condition is equivalent to the condition
$$ {\cal L}_{\cal J} (W({\cal J}^* \omega^2)) - {\cal J}^* (W({\cal L}_{\cal J} \omega^2)) + W({\cal L}_{\cal J} \omega^2) - {\cal L}_{\cal J} (W(\omega^2))=0 , $$

where we define, for any 2-form $\omega^2 \in \Lambda^2 (M^{16})$ and any complex structure ${\cal J} \in {\cal T}_1 (M^{16})$, the 2-forms ${\cal J}^* (\omega^2)$ and ${\cal L}_{\cal J} (\omega^2)$ by the formulas:
\begin{eqnarray*}
({\cal J}^* \omega^2)(X,Y) &=& \omega^2 ({\cal J} X, {\cal J}Y)\\
({\cal L}_{\cal J} \omega^2)(X,Y) &=& \omega^2 ({\cal J} X,Y) + \omega^2 (X, {\cal J} Y) . 
\end{eqnarray*}

Since ${\cal J}$ is an antisymmetric complex structure, we have the following relations on 2-forms:
$$ {\cal J}^* ({\cal L}_{\cal J} \omega^2)= {\cal L}_{\cal J} ({\cal J}^* \omega^2)= - {\cal L}_{\cal J} \omega^2 \quad , \quad  ({\cal J}^* )^2 \omega^2 = \omega^2 . $$

The $\sneun$-structure yields a splitting of the bundle of 2-forms:
$$ \Lambda^2 (\ms)= \Lambda^2 (V^9) \oplus \Lambda^3 (V^9) . $$

It turns out that, for ${\cal J} \in {\cal T}_1 (\ms)$, the operators ${\cal J}^*$ and ${\cal L}_{\cal J}$ preserve this splitting. Let us thus write the endomorphism $W: \Lambda^2 (\ms) \to \Lambda^2 (\ms)$ as a $(2 \times 2)$-matrix
$$ W= \left( \begin{array}{cc} W_{22} & W_{32} \\ W_{23} & W_{33} \end{array} \right)  , $$

where $W_{\alpha \beta} : \Lambda^{\alpha} (V^9) \to \Lambda^{\beta} (V^9) \quad (\alpha, \beta =2,3)$ is a bundle morphism. Then the integrability condition splits into 4 conditions, too:
$$ {\cal L}_{\cal J} (W_{\alpha \beta} ({\cal J}^* \omega^2)) - {\cal J}^* (W_{\alpha \beta} ({\cal L}_{\cal J} \omega^2)) + W_{\alpha \beta} ({\cal L}_{\cal J} \omega^2) - {\cal L}_{\cal J} (W_{\alpha \beta} (\omega^2))=0 . $$

We can analyze these conditions in the usual way (see for example \cite{OBR}) using representation theory. However, since $W$ does not satisfy the first Bianchi identity,  the discussion becomes more complicated. We will not provide this discussion in details, but let us  investigate the condition for $W_{22}$ for example. The complex structure ${\cal J} \in {\cal C} \cap \Lambda^2 ({\Bbb R}^9) = {\cal T}_1$ is an element of the Lie algebra $\mbox{\graf spin}(9)$ as well as of the group $\sneun$. In case the two form $\omega^2$ belongs to $\Lambda^2 (V^9)$, we have
$$ {\cal L}_{\cal J} (\omega^2) = - [ {\cal J} , \omega^2 ] \quad , \quad {\cal J}^* (\omega^2)= Ad ({\cal J} )(\omega^2) . $$

Therefore $W_{22}: \mbox{\graf spin}(9) \to \mbox{\graf spin}(9)$ is a linear map satisfying the condition 
$$ [€{\cal J}, W_{22} (Ad({\cal J})\omega^2)] - Ad ({\cal J})(W_{22} ([{\cal J}, \omega^2 ])) + W_{22}([{\cal J}, \omega^2]) - [{\cal J} , W_{22}(\omega^2)]=0 $$

for any $\omega^2 \in \Lambda^2 (V^9)$ and any ${\cal J} \in {\cal C} \cap \mbox{\graf spin}(9)$. The adjoint action of the complex structure ${\cal J}$ on a 2-form $\omega^2$ can be expressed by the commutator
$$ Ad ({\cal J}) \omega^2 = \omega^2 + \frac{1}{2} [{\cal J} , [{\cal J}, \omega^2 ]] . $$

Indeed, we have
\begin{eqnarray*}
[{\cal J}, [{\cal J}, \omega^2 ]] &=& {\cal J}({\cal J} \omega^2 - \omega^2 {\cal J}) - ({\cal J}\omega^2 - \omega^2 {\cal J}) {\cal J} =\\
&=& - 2 \omega^2 + 2 {\cal J} \omega^2 {\cal J}^* = - 2 \omega^2 + Ad ({\cal J}) \omega^2 . 
\end{eqnarray*}

Using this formula, the integrability condition for $W_{22}$ becomes much simpler:
$$ [ {\cal J}, W_{22} ([{\cal J}, [{\cal J}, \omega^2]])] =[ {\cal J} , [{\cal J}, W_{22} ([{\cal J}, \omega^2 ])]] . $$

Fix a complex structure ${\cal J} \in {\cal T}_1$ and decompose the Lie algebra $\mbox{\graf spin}(9)$ into the stabilizer of ${\cal J}$
$$ \h_{{\cal J}} = \Big\{ \omega^2 \in \mbox{\graf spin}(9): \quad [{\cal J} , \omega^2 ]=0 \Big\} $$

and its orthogonal complement  $\h_{{\cal J}}^{\perp}$. This subspace $\h_{\cal J}^{\perp}$ is the tangent space $T_{\cal J} ({\cal T}_1)$ of the quadric ${\cal T}_1$ at the point ${\cal J} \in {\cal T}_1$. Moreover, $\frac{1}{2} ad ({\cal J})$ defines the complex structure of ${\cal T}_1$. In fact, for any 2-form $\omega^2 \in \h_{\cal J}^{\perp}$ the formula
$$ [ {\cal J} , [{\cal J}, \omega^2 ]] = - 4 \omega^2 $$

holds. The latter equation defines $\h_{\cal J}^{\perp}$ as a subspace of $\mbox{\graf spin}(9)$:
$$ \h_{\cal J}^{\perp} = \Big\{ \omega^2 \in \mbox{\graf spin}(9): \quad [{\cal J}, [{\cal J}, \omega^2]]= - 4 \omega^2 \Big\} . $$

The quadric ${\cal T}_1$ is a symmetric space. Therefore, in the decomposition
$$ \mbox{\graf spin}(9) = \h_{\cal J} \oplus \h_{\cal J}^{\perp} $$

the commutator relations
$$ [ \h_{\cal J}, \h_{\cal J}] \subset \h_{\cal J}  \quad , \quad [\h_{\cal J} , \h_{\cal J}^{\perp}] \subset \h_{\cal J}^{\perp} \quad , \quad [h_{\cal J}^{\perp} , \h_{\cal J}^{\perp} ] \subset h_{\cal J} $$

hold. We remark that for any real number $c \in {\Bbb R}^1$, any 2-form $\eta^2 \in \Lambda^2 ({\Bbb R}^9)$ and any 5-form $\mu^5 \in \Lambda^5 ({\Bbb R}^9)$ the endomorphism $W_{22}: \Lambda^2 ({\Bbb R}^9) \to \Lambda^2 ({\Bbb R}^9)$ defined by the formula
$$ W_{22} (\omega^2)= c \cdot \omega^2 + [\eta^2 , \omega^2 ] + \star (\mu^5 \wedge \omega^2) $$

satisfies the integrability condition. Indeed, the $\eta^2$-term can be handled using the Jacobi identity. Then we obtain $(\mu^5 =0)$
$$ \Big[{\cal J}, W_{22} ([{\cal J}, [{\cal J}, \omega^2 ]]) - [{\cal J}, W_{22} ([ {\cal J}, \omega^2 ])]\Big] = \Big[{\cal J} , [[{\cal J} , \omega^2 ], [{\cal J}, \eta^2 ]]\Big] . $$

The elements $[{\cal J}, \omega^2]$ and $[{\cal J}, \eta^2 ]$ belong to $h_{\cal J}^{\perp}$ and, consequently, $[[{\cal J} , \omega^2], [{\cal J} , \eta^2]]$ is an element of the Lie algebra $\h_{\cal J}$. The $\mu^5$-term in the formula of $W_{22}$ satisfies the integrability condition, too. This is a consequence of the algebraic relations of the endomorphisms $I_{\alpha} \, \, (1 \le \alpha \le 9)$ and can be checked by a direct calculation. Altogether we obtain a family of endomorphisms $W_{22}$ satisfying the integrability condition and depending on
$$ \dim \Lambda^2 ({\Bbb R}^9) + \dim \Lambda^5 ({\Bbb R}^9) + 1 = 163 $$

parameters. A representation-theoretic argument shows that we derived the general formula for the endomorphism $W_{22}$. Let us sketch the argument. The tensor product $\Lambda^2 ({\Bbb R}^9) \otimes \Lambda^2 ({\Bbb R}^9)$ decomposes into six irreducible $SO(9)$-representations:
\begin{eqnarray*}
\Lambda^2 ({\Bbb R}^9) \otimes \Lambda^2 ({\Bbb R}^9) &=& \langle 1 \rangle \oplus S_0^2 (\Lambda^2 ({\Bbb R}^9)) \oplus \Lambda^2 (\Lambda^2 ({\Bbb R}^9) ) =\\
&=& \langle 1 \rangle  \oplus \{ \langle 44 \rangle \oplus \langle 126 \rangle \oplus \langle 495 \rangle \} \oplus \{ \langle 36 \rangle \oplus \langle 594 \rangle \} 
\end{eqnarray*}

where, for example, $\langle 44 \rangle$ denotes the unique irreducible $SO(9)$-representation of dimension 44. We have already proved that the terms $\langle 1 \rangle = {\Bbb R}^1, \langle 36 \rangle = \Lambda^2 ({\Bbb R}^9)$ and $\langle 126 \rangle = \Lambda^5 ({\Bbb R}^9)$ may occur in the decomposition of $W_{22}$. Therefore, we have to exclude the representations $\langle 44 \rangle, \langle 495 \rangle$ and $\langle 594 \rangle $. A traceless symmetric endomorphism $L: {\Bbb R}^9 \to {\Bbb R}^9$ induces by
$$ {\cal L}_L \omega^2 (\alpha , \beta)= \omega^2 (\mbox{L} \alpha, \beta) + \omega^2 (\alpha, \mbox{L} \beta)$$

a linear map ${\cal L}_L: \Lambda^2 ({\Bbb R}^9) \to \Lambda^2 ({\Bbb R}^9)$ and this formula realizes the irreducible representation $\langle 44 \rangle =S_0^2 ({\Bbb R}^9)$ in $\Lambda^2 ({\Bbb R}^9) \otimes \Lambda^2 ({\Bbb R}^9)$. Consider, in particular, the endomorphism  $\mbox{L}_0 \in S_0^2 ({\Bbb R}^9)$ defined by the formulas
$$ L_0 (I_1)=I_2 \quad , \quad L_0 (I_2)= I_1 \quad , \quad L_0 (I_{\alpha})=0 \, \, \, \, \mbox{for} \, \, \, 3 \le \alpha \le 9 $$

as well as the complex structure ${\cal J}_0 =I_1 I_2$ and the 2-form $\omega^2_0 =I_1 \wedge I_3$. A direct computation yields the result
$$ [{\cal J}_0 , {\cal L}_{\mbox{L}_0} ([{\cal J}_0, [{\cal J}_0 , \omega_0^2 ]])]= - 8 I_1 I_3  \quad , \quad [{\cal J}_0 , [{\cal J}_0, {\cal L}_{\mbox{L}_0} ([{\cal J}_0 , \omega^2_0 ])]]= 8 I_1 I_3 . $$

Consequently, the endomorphism $W_{22}$ is orthogonal to the representation $\langle 44 \rangle$. Let us summarize  the previous discussion:\\

{\bf Theorem 2:} {\it An endomorphism $W_{22}: \Lambda^2 ({\Bbb R}^9) \to \Lambda^2 ({\Bbb R}^9)$ satisfies the integrability condition if and only if there exists a constant $c \in {\Bbb R}^1$, a 2-form $\eta^2 \in \Lambda^2 ({\Bbb R}^9)$ and a 5-form $\mu^5 \in \Lambda^5 ({\Bbb R}^9)$ such that 
$$ W_{22} (\omega^2) = c \cdot \omega^2 + [\eta^2 , \omega^2 ] + \star (\mu^5 \wedge \omega^2)$$

holds.}\\

In a similar way one can discuss the possible type of the endomorphisms $W_{23}, W_{32}$ and $W_{33}$.\\

The twistor space ${\cal T}_1 (\ms)$ is a fibration over $\ms$ and the fibres are complex submani\-folds analytically isomorphic to the quadric ${\cal T}_1 =Q$ in ${\Bbb P}^8 ({\Bbb C})$. The group $Iso (\ms ; {\cal R})$ of all isometries of $\ms$ preserving the $\sneun$-structure ${\cal R} \subset {\cal F} (\ms)$ acts on the twistor space as a group of holomorphic transformations.\\

{\bf Example 1:} (The twistor space of the Cayley plane)  \\
The Cayley plane $F_4 /\sneun$ is a symmetric space and the Riemannian connection reduces to the $\sneun$-structure. The twistor space
$$ {\cal T}_1 (F_4 / \sneun )=F_4 / (Spin(2) \, \mbox{x}_{{\Bbb Z}_2} Spin(7)) $$

is a 15-dimensional complex manifold. The exceptional group $F_4$ acts transitively on the twistor space as a group of holomorphic transformations. Consider the torus $T_1= Spin(2)$. Its centralizers in $\sneun$ and in $F_4$ coincide:
$$ C_{\sneun} (T_1)= C_{F_4} (T_1) = Spin(2) \, \mbox{x}_{{\Bbb Z}_2} Spin(7) . $$

Consequently, the twistor space ${\cal T}_1 (F_4 / \sneun)$ of the Cayley plane is a generalized flag manifold and, henceforth, a projective variety (see \cite{wallach}).\\

{\bf Example 2:} (The twistor space of $S^1 \times S^{15}$)\\
$S^1 \times S^{15}$ admits a nearly parallel $\sneun$-structure with the parallel vector field $\Gamma = \frac{\partial}{\partial \varphi}$ tangent to $S^1$. Therefore, the tensor field $W$ coincides with the curvature tensor ${\cal R}$ of the space form $S^{15}$. Using the well-known formula
$$ {\cal R} (X,Y) Z= \langle Y, Z \rangle X - \langle X,Z \rangle Y $$

we immediately see that the integrability condition of Theorem 1 holds, i.e. \linebreak ${\cal T}_1 (S^1 \times S^{15})$ is a complex manifold. Since the first Betti number is odd, ${\cal T}_1 (S^1 \times S^{15})$ does not admit a K\"ahler metric. In particular, this complex manifold is not algebraic projective.\\

{\bf Example 3:}  (The twistor space of ${\Bbb R}^{16}$)\\
Since ${\Bbb R}^{16}$ is flat, the parallel displacement along lines through a fixed point $0 \in {\Bbb R}^{16}$ defines a holomorphic projection of the twistor space ${\cal T}_1 ({\Bbb R}^{16})$ onto one fibre $Q$. Consequently, ${\cal T}_1 ({\Bbb R}^{16})$ is analytically isomorphic to an 8-dimensional holomorphic vector bundle $N$ over $Q$. Let us describe the bundle $N$. Consider the trivial real vector bundle
$$ N=Q \times {\Bbb R}^{16} $$

and introduce a complex structure on $N$ by the formula
$$ i \cdot ({\cal J} , \mbox{x}):= ({\cal J}, {\cal J} \mbox{x}) $$

for ${\cal J} \in Q= {\cal T}_1$ and $\mbox{x} \in {\Bbb R}^{16}$. The group $\sneun$ acts on the bundle $N$ by
$$ g \cdot ({\cal J} , \mbox{x}):=(g {\cal J} g^{-1} , g \mbox{x} )$$

and, therefore, $N$ is a homogeneous vector bundle over the space
$$  {\cal T}_1 = \sneun / Spin(2) \mbox{x}_{{\Bbb Z}_2} Spin(7)  \quad , \quad N= \sneun \mbox{x}_{Spin(2) \mbox{x}_{{\Bbb Z}_2} Spin(7)} {\Bbb R}^{16} . $$

As usual, the bundle $N$ is the associated bundle to the representation of the group $Spin(2) \mbox{x}_{{\Bbb Z}_2} Spin(7)= Spin(9) \cap U(8)$ in $U(8)$. We decompose the Lie algebra $\mbox{\graf spin}(9)$ 
$$ \mbox{\graf spin}(9) = (\mbox{\graf spin}(2) \oplus \mbox{\graf spin}(7)) \oplus \mbox{\graf l} $$

where the linear space $\mbox{\graf l}$ consist of all elements
$$ \sum\limits^9_{\alpha =3} \mbox{x}_{1 \alpha} I_1 I_{\alpha} + \sum\limits^9_{\beta =3} \mbox{x}_{2 \beta} I_2 I_{\beta} . $$

We can compute the $\mbox{\graf u}(8)$-valued curvature form $\Omega^N$ of the bundle $N$ using the canonical connection of the symmetric space $Q$. This formula is quite a  complicated one. However, the trace of $\Omega^N$ has a simple form:
$$ Tr (\Omega^N)= - 16i (d\mbox{x}_{13} \wedge d\mbox{x}_{23} + \ldots + d\mbox{x}_{19} \wedge d\mbox{x}_{29}) . $$

Since the first Chern class $c_1 (N)$ is represented by the form $c_1 (N)= - \frac{1}{2 \pi i} Tr (\Omega^N)$, we obtain a formula for $c_1 (N)$:
$$ c_1 (N)= \frac{16}{2 \pi} (d\mbox{x}_{13} \wedge d \, \mbox{x}_{23} + \ldots + d \, \mbox{x}_{19} \wedge d \, \mbox{x}_{29} ) . $$

On the other hand, $Q$ is a K\"ahler-Einstein manifold and therefore the first Chern class $c_1 (Q)$ is proportional to the fundamental form of the K\"ahler structure. In this way we obtain the formula
$$ c_1 (Q) = \frac{14}{2 \pi} (d\mbox{x}_{13} \wedge d\mbox{x}_{23} + \ldots + d\mbox{x}_{19} \wedge d\mbox{x}_{29} ) $$

and, finally, we have the relation
$$ c_1 (N)= \frac{8}{7} \, c_1 (Q) . $$

Let $\mbox{x} \in H^2 ({\Bbb P}^8 ({\Bbb C}); {\Bbb Z})$ be the generator of the second cohomology group of ${\Bbb P}^8 ({\Bbb C})$. The Chern classes $c_i (Q)$ of the quadric $Q \subset {\Bbb P}^8 ({\Bbb C})$ can be described by the powers of the element $\mbox{x}$:
$$c_1 (Q)= 7 \,  \mbox{x} \quad \hspace{0.2cm} , \quad \hspace{0.1cm} c_2 (Q)= 22 \, \mbox{x}^2 \quad , \quad c_3 (Q)= 40 \, \mbox{x}^3 \hspace{2.6cm} \mbox{}$$

\mbox{} \hspace{1.4cm} $c_4 (Q)= 46 \, \mbox{x}^4  \hspace{0.2cm} , \quad \, \,  c_5 (Q)= 34 \, \mbox{x}^5 \quad , \quad c_6 (Q)= 16 \,  \mbox{x}^6 \quad , \quad c_7 (Q)= 4 \, \mbox{x}^7 $\\

(see \cite{7x}).  Consequently, we obtain the formula
$$ c_1 (N)= 8 \, \mbox{x} . $$

The real vector bundle $N$ is trivial and its first Pontrjagin class $p_1 (N)$ vanishes:
$$ 0= p_1 (N) = 2 \, c_2 (N) - c_1^2 (N) . $$

Therefore, the second Chern class $c_2 (N)$ is proportional to $c_1^2 (N)$:
$$ c_2 (N) = \frac{1}{2} \, c_1^2 (N) = 32 \, \mbox{x}^2 . $$

We describe the space ${\cal H}^0 (N)$ of all holomorphic sections of the bundle $N$ over $Q$. Any fixed vector $y \in {\Bbb R}^{16} = \dneun$ defines a section $S_y \in \Gamma (Q;N)$ by the map $S_y : \sneun \to {\Bbb R}^{16} = \dneun$,
$$ S_y (A)= A^{-1} (y) . $$

$S_y$ is a holomorphic section and the $\sneun$-action on ${\cal H}^0 (N)$ coincides under this identification with the $\sneun$-action on ${\Bbb R}^{16}$. In particular, we have computed the dimension of ${\cal H}^0 (N)$:
$$ \dim_{\Bbb C}  {\cal H}^0 (N) = 16 , $$

and ${\cal H}^0 (N)$ is the unique irreducible $\sneun$-representation of dimension 16.
\mbox{} \hfill \rule{3mm}{3mm}\\

The previous discussion describes not only the twistor space ${\cal T}_1 ({\Bbb R}^{16})$ of the flat space, but also the normal bundle to any fibre inside an arbitrary twistor space ${\cal T}_1 (\ms)$.\\

{\bf Theorem 3:} {\it Let $\ms$ be a 16-dimensional Riemannian manifold with a nearly parallel $\sneun$-structure and suppose that the twistor space ${\cal T}_1 (\ms)$ is analytic. The normal bundle $N$ to any fibre $Q \subset {\cal T}_1 (\ms)$ is an 8-dimensional holomorphic vector bundle over the quadric $Q$ with the following properties:}
\begin{enumerate}
\item $\displaystyle c_1 (N)= \frac{8}{7} c_1 (Q)= 8 \, \mbox{x} \quad , \quad c_2 (N)= \frac{1}{2} c_1^2 (N) = 32 \, \mbox{x}^2 . $
\item $\displaystyle \dim_{\Bbb C} {\cal H}^0 (N)= 16$.
\end{enumerate}
\mbox{} \hfill \rule{3mm}{3mm}\\

The fibre ${\cal T}_1$ of the twistor space ${\cal T}_1 (\ms)$ admits an antiholomorphic involution $\tau$. Using the different models for ${\cal T}_1$ we can define $\tau$ in different ways. The involution $\tau : {\cal T}_1 \to {\cal T}_1$ is given by 
$$ \tau ({\cal J})= - {\cal J} . $$

In case we identify ${\cal T}_1$ with the quadric
$$ Q= \Big\{ [ z_0 : \ldots : z_8 ] \in {\Bbb P}^8 ({\Bbb C}) : z_0^2 + \ldots + z_8^2 =0 \Big\} , $$

the involution $\tau$ is the conjugation, $ \tau [z_0 : z_1 : \ldots : z_8 ]=[ \bar{z}_0 : \bar{z}_1 : \ldots : \bar{z}_8 ]$. Finally, in case ${\cal T}_1 = \sneun / Spin(2) \mbox{x}_{{\Bbb Z}_2} Spin(7)$ is considered as the Gra\ss mann manifold $G_{2,9}$ of all oriented 2-planes $\pi^2$ in ${\Bbb R}^9$,  the involution $\tau$ reverses the orientation, $ \tau (\pi^2)= - \pi^2$. Since $\tau$ commutes with the $\sneun$-action on ${\cal T}_1$, it defines an involution
$$ \tau : {\cal T}_1 (\ms) \to {\cal T}_1 (\ms) $$

on any twistor space. The map $\tau$ is an antiholomorphic involution without fixed points. Moreover, $\tau$ acts on the space ${\cal H}^0 (N)$ of all holomorphic sections of the normal bundle to a fibre $Q \subset {\cal T}_1 (\ms)$ and the space of real sections can be identified with the tangent space at the base point to $\ms$.\\

Summarizing, the twistor space ${\cal T}_1 (\ms)$ of a 16-dimensional manifold with a (nearly parallel) $\sneun$-structure has the same structure as the twistor space of an oriented 4-dimensional Riemannian manifold. The difference is the more complicated topology of the fibre. \\

{\bf Remark:} We defined a twistor space ${\cal T}_1 (\ms)$ by using the space ${\cal T}_1 {\cal C} \cap \Lambda^2 ({\Bbb R}^9)$ of all complex structures that are given by a two-form in $\Lambda^2 ({\Bbb R}^9)$. There is a second possibility. Consider the space
$$ {\cal T}_2 = {\cal C} \cap \Lambda^3 ({\Bbb R}^9) = \Big\{ {\cal J} = \sum\limits_{1 \le \alpha < \beta < \gamma \le 9} y_{\alpha \beta \gamma} \, \, \, I_{\alpha} I_{\beta} I_{\gamma} : \, \, \, {\cal J}^2 = - \mbox{Id}  \Big\} $$

of all complex structures on ${\Bbb R}^{16}$ defined by a 3-form in $\Lambda^3 ({\Bbb R}^9)$. Then the group $\sneun$ does not act transitively on ${\cal T}_2$. For example, $I_1 I_2 I_3$ and $\frac{1}{\sqrt{3}} (I_1 I_2 I_3 + I_4 I_5 I_6 + I_7 I_8 I_9)$ are two elements in ${\cal T}_2$ with different isotropy groups with respect to the $\sneun$-action. The complete orbit structure of ${\cal T}_2$ is very difficult and related to the classification of all normal forms of 3-forms in 9 variables. For the  case of the group $SL(9)$ acting on ${\Bbb R}^9$, the orbit structure of $\Lambda^3 ({\Bbb R}^9)$ was described in the paper \cite{vinberg}. Anyway, consider a $\sneun$-orbit ${\cal O}^* \subset {\cal T}_2$. Then we can introduce the twistor space
$$ {\cal T}^{{\cal O}^*}_2 (\ms) = {\cal R} \times_{\sneun} {\cal O}^* $$

and a given $\sneun$-invariant geometric structure on ${\cal O}^*$ induces a similar structure on the twistor space ${\cal T}^{{\cal O}^*}_2$.\\

\vskip 1.0 cm  

\section{References}  

\vspace{-0.5cm}
\begin{Literatur}{MMM}
\bibitem {agricola} I. Agricola, On the decomposition of antisymmetric tensor products of some spin  representations, to appear.
\bibitem {11x} D.V. Alekseevski, Riemannian spaces with exceptional holonomy groups, Func. Anal. Prilozh. 2 (1968), 1-10.
\bibitem {BFGK} H. Baum, Th. Friedrich, R. Grunewald, I. Kath, Twistors and Killing spinors on Riemannian manifolds, Teubner-Verlag Stuttgart-Leipzig 1991.
\bibitem {BBO} L. Berard-Bergery, T. Ochiai, On some generalizations of the construction of twistor spaces, Durham-LMS Symposium, Global Riemannian Geometry (ed. by T.J. Willmore and N. Hitchin), New York 1984, 52-59.
\bibitem {6x} A. Borel, F. Hirzebruch, Characteristic classes and homogeneous spaces I, Amer. Journ. Math. 80 (1958), 458-538.
\bibitem {x} R.B. Brown, A. Gray, Riemannian manifolds with holonomy group Spin(9), Differential Geometry in honor of K. Yano, Kinokuniya, Tokyo 1972, 41-59.
\bibitem {4x} E. Cartan, Le\c{c}ons sur la th\'eorie de spineurs, Hermann, Paris 1938.
\bibitem {F2} M. Fernandez, A classification of Riemannian manifolds with structure group Spin(7), Ann. Mat. Pura Appl. 143 (1986), 101-122.
\bibitem {FG1} M. Fernandez, A. Gray, Riemannian manifolds with structure group $G_2$, Ann. Mat. Pura Appl. 132 (1982), 19-45.
\bibitem {3x} A.T. Fomenko, Variational principles in topology,  Kluwer Academic Publishers 1990.
\bibitem {xx} Th. Friedrich, Dirac-Operatoren in der Riemannschen Geometrie, Vieweg-Verlag, Braunschweig/Wiesbaden 1997. 
\bibitem {FKMS} Th. Friedrich, I. Kath, A. Moroianu, U. Semmelmann, On nearly parallel $G_2$-structures, J. Geom. Phys. 27 (1998), 155-177.
\bibitem {Gray1} A. Gray, Nearly K\"ahler manifolds, J. Diff. Geom. 4 (1970), 283-309.
\bibitem {Gray2} A. Gray, Weak holonomy groups, Math. Zeitschrift 123 (1971), 290-300.
\bibitem {Gray3} A. Gray, The structure of nearly K\"ahler manifolds, Math. Ann. 223 (1976), 233-248.
\bibitem {5x} A. Gray, P. Green, Sphere transitive structures and the triality automorphism, Pac. J. Math. 34 (1970), 83-96.
\bibitem {Grune} R. Grunewald, Six-dimensional Riemannian manifolds with a real Killing spinor, Ann. Glob. Anal. Geom. 8 (1990), 43-59.
\bibitem {7x} F. Hirzebruch, Topological methods in algebraic geometry, Springer-Verlag 1978.
\bibitem {huse} D. Husemoller, Fibre bundles, New York 1966.
\bibitem {9x} A.T. Lundell, Concise tables of James numbers  and some homotopy of classical Lie groups and associated homogeneous spaces, Algebraic topology, Proc. Conf. S. Felin de Guixols/ Spain 1990, Lect. Notes Math. 1509 (1992), 250-272.
\bibitem {March} S. Marchiafava, Characterization of Riemannian manifolds with weak holonomy group $G_2$, Math. Zeitschrift 178 (1981), 158-162.
\bibitem {OBR} N.R. O'Brian, J.H. Rawnsley, Twistor spaces, Ann. Glob. Anal. Geom. 3 (1985), 29-58.
\bibitem {grady} K. O'Grady, The Hodge structure of the intersection of three quadrics in an odd dimensional projective space, Math. Ann. 273 (1985), 277-285.
\bibitem {8x} I.R. Shafarevich, Basic algebraic geometry, vol. 2, Springer-Verlag 1994.
\bibitem {10x} M.J. Slupinski, A Hodge type decomposition for spinor valued forms, Ann. scient. \'Ec. Norm. Sup., $4^e$ s\'erie, t. 29 (1996), 23-48.
\bibitem {vinberg} E.B. Vinberg, A.G. Elashvili, Classification of trivectors of a 9-dimensional space, Sel. Math. Sov. vol. 7 (1988), 63-98; originally published in Trudy Sem. Vektor. Tenzor. Anal. 18 (1978), 197-233.
\bibitem {wallach} N.R. Wallach, Harmonic analysis on homogeneous spaces, New York 1973.
\end{Literatur}

\end{document}